\documentclass[11pt,a4paper]{amsart}
\usepackage{amscd,amssymb,amsopn,amsmath,amsthm,graphics,amsfonts,enumerate,verbatim,calc}
\usepackage[dvips]{graphicx}

\usepackage{amssymb,amsmath}
\usepackage{mathpazo}

\headsep=1cm \numberwithin{equation}{section}
\newtheorem{theorem}{Theorem}[section]
\newtheorem{lemma}[theorem]{Lemma}
\newtheorem{proposition}[theorem]{Proposition}
\newtheorem{corollary}[theorem]{Corollary}
\newtheorem{claim}[theorem]{Claim}

\theoremstyle{definition}
\newtheorem{remark}[theorem]{Remark}
\newtheorem{discussion}[theorem]{Discussion}
\newtheorem{definition}[theorem]{Definition}
\theoremstyle{remark}
\newtheorem{example}[theorem]{Example}
\newtheorem{question}[theorem]{Question}

\newtheorem{acknowledgement}{Acknowledgement}

\newtheorem{maintheorema}{Main Theorem A}

\newtheorem{maintheoremb}{Main Theorem B}

\newcommand{\Ass}{\operatorname{Ass}}

\newcommand{\Spec}{\operatorname{Spec}}

\newcommand{\cl}{\operatorname{rc}}
\newcommand{\Ht}{\operatorname{ht}}

\newcommand{\Ext}{\operatorname{Ext}}
\newcommand{\Supp}{\operatorname{Supp}}

\newcommand{\Hom}{\operatorname{Hom}}

\newcommand{\Ann}{\operatorname{Ann}}
\newcommand{\Proj}{\operatorname{Proj}}

\newcommand{\Char}{\operatorname{char}}
\newcommand{\Fitt}{\operatorname{Fitt}}
\newcommand{\red}{\operatorname{red}}

\newcommand{\depth}{\operatorname{depth}}

\newcommand{\coker}{\operatorname{coker}}

\newcommand{\Min}{\operatorname{Min}}
\newcommand{\rank}{\operatorname{rank}}
\newcommand{\Sp}{\operatorname{Sp}}
\newcommand{\ur}{\operatorname{ur}}

\newcommand{\fm}{\frak{m}}
\newcommand{\fp}{\frak{p}}
\newcommand{\fq}{\frak{q}}
\newcommand{\fa}{\frak{a}}

\newcommand{\fn}{\frak{n}}

\begin{document}
\title[Bertini theorem for normality on local rings in mixed characteristic]
{Bertini theorem for normality on local rings in mixed characteristic (applications to characteristic ideals)}

\author[T.Ochiai]{Tadashi Ochiai}
\address{Graduate school of Science
Osaka University 1-1 Machikaneyama Toyonaka 
Osaka 560-0043 Japan}
\email{ochiai@math.sci.osaka-u.ac.jp}
\thanks{The first-named author is partially supported by Grant-in-Aid for Challenging Exploratory Research (24654004) and Sumitomo Foundation.}

\author[K. Shimomoto]{Kazuma Shimomoto}
\address{Department of Mathematics, School of Science and Technology Meiji Univesity 1-1-1 Higashimita Tama-Ku Kawasaki 214-8571 Japan}
\email{shimomotokazuma@gmail.com}
\thanks{The second-named author is partially supported by Grant-in-Aid for Young Scientists (B) \# 25800028}




\begin{abstract}
In this article, we prove a strong version of local Bertini theorem for normality on local rings in mixed characteristic. The main result asserts that a generic hyperplane section of a normal, Cohen-Macaulay, and complete local domain of dimension at least 3 is normal. Applications include the study of characteristic ideals attached to torsion modules over normal domains, which is fundamental in the study of Euler system theory, Iwasawa's main conjectures, and the deformation theory of Galois representations.
\end{abstract}

\maketitle

\section{Introduction}

The classical Bertini theorem says that a generic hyperplane section of a smooth complex projective variety is smooth. 
We would like to call by a "local Bertini theorem" the following problem for a local Noetherian ring $(R,\fm,\mathbf{k})$: 

\begin{quotation}
Let $\mathbf{P}$ be a ring-theoretic property (e.g. regular, reduced, normal, seminormal and so on). Then if $R$ is $\mathbf{P}$ and $x \in \fm$ is a non-zero divisor, then is it true that $R/xR$ is $\mathbf{P}$ for a generic choice of $x$? 
\end{quotation}

We refer the reader to Definition \ref{definition:Bertini_property} for a more precise formulation of the problem. A local Bertini theorem (in a slightly weak form) was first raised by Grothendieck (\cite{Gro}, Expos\'e XIII, Conjecture 2.6) and was proved by Flenner \cite{Fle} and Trivedi \cite{Triv}. We use the following notation. For an ideal $I \subseteq R$, we denote by $U(I)$ the set of all primes of $R$ which do not contain $I$, and by $V(I)$ the complement of $U(I)$ in $\Spec R$. For a graded ring $R$ with a graded ideal $I$, we use the notation $V^+(I)$ and $U^+(I)$ for a closed subset with its complement in $\Proj R$. Before stating our main theorems, let us recall the following result from \cite{Fle}:

\begin{theorem}[Flenner-Trivedi]
\label{theorem1}
Let $(R,\fm)$ be a local Noetherian ring and let $I \subseteq \fm$ be an ideal. Assume that $Q$ is a finite subset of $U(I)$. Then there exists an element $x \in I$ such that:

\begin{enumerate}
\item[$\mathrm{(1)}$]
$x \notin \fp^{(2)}$ for all $\fp \in U(I)$;

\item[$\mathrm{(2)}$]
$x \notin \fp$ for all $\fp \in Q$;
\end{enumerate}
for $\fp^{(n)}:=\fp^n R_{\fp} \cap R$, the $n$-th symbolic power ideal of $\fp$.
\end{theorem}

We make a remark on the second symbolic power of ideals. Let $(R,\fm)$ be a local ring, let $x \in \fp$ be a non-zero divisor and let $x \notin \fp^{(2)}$ for a prime ideal $\fp \subseteq R$. Then $R_{\fp}$ is regular if and only if $R_{\fp}/xR_{\fp}$ is regular. Many ring-theoretic properties such as regular, normal, reduced can be verified at the localization $R_{\fp}$, which is the reason why we require $x \notin \fp^{(2)}$ (but not merely $x \notin \fp^2$), which is equivalent to the condition: $x \notin \fp^2 R_{\fp} \cap R$. A strong version of local Bertini theorem similar to our main theorem below was already proved for local rings containing a field in \cite{Fle}. To extend it to the mixed characteristic case, we need to introduce some new ideas. Theorem \ref{theorem1} has the following implication. If $(R,\fm)$ is a complete local normal domain of $\depth R \ge 3$, then there is a non-zero element $x \in \fm$ for which $R/xR$ is normal (\cite{Fle}, Korollar 3.4). However, this does not suffice for our purpose, because we would like to give a family of infinitely many such specializations $R/xR$ which are parametrized by a certain explicit $p$-adic space.

Here is our notation which will be used throughout the paper. 
We denote by $\mathbb{P}^{n} (S)$ (resp. $\mathbb{A}^n(S)$) the set of $S$-rational points of the $n$-dimensional projective space (resp. affine space) over an integral domain $S$. In case when $S$ is a field or a discrete valuation ring, these spaces come equipped with some topology (see discussions in \S~\ref{topology}). For a complete discrete valuation ring $(A,\pi_A,\mathbf{k})$, we have a \textit{specialization map} $\Sp_A:\mathbb{P}^{n}(A) \to \mathbb{P}^{n}(\mathbf{k})$ (cf. Definition \ref{definition:basic_p_adic}). 

In \S~\ref{topology} and \S~\ref{basicelement}, we make discussions on basic tools and prove some requisite results, including the specialization map and basic elements, due to Swan.

In \S~\ref{sectionmain}, we prove the first main theorem (see Theorem \ref{Bertini} together with Theorem \ref{theorem2}) in this article:

\begin{maintheorema}[Local Bertini Theorem]
Let $(R,\fm,\mathbf{k})$ be a complete local domain of mixed characteristic $p>0$ and suppose the following conditions:

\begin{enumerate} 
\item[$\mathrm{(1)}$]
let $A \to R$ be a coefficient ring map for a complete discrete valuation ring $(A,\pi_{A})$; 

\item[$\mathrm{(2)}$]
let $x_0,x_1,\ldots,x_d$ be a fixed set of minimal generators of $\fm$;

\item[$\mathrm{(3)}$]
$R$ is normal, of $\depth R \ge 3$, and the residue field $\mathbf{k}$ is infinite.
\end{enumerate}
Then there exists a Zariski dense open subset $U \subseteq \mathbb{P}^{d}(\mathbf{k})$ satisfying the following properties. For any $a=(a_{0}:\cdots :a_{d}) \in \Sp_A^{-1}(U)$, the quotient $R/\mathbf{x}_{\widetilde{a}}R$ is a normal domain of mixed characteristic $p>0$ (cf. Definition \ref{definition:Bertini_property} for the notation $\mathbf{x}_{\widetilde{a}}R$). 
\end{maintheorema}

By the minimality condition of $x_0,x_1,\ldots,x_d$, it is immediate to see that $\mathbf{x}_{\widetilde{a}} \ne 0$, that is, $R/\mathbf{x}_{\widetilde{a}}R$ is a non-trivial quotient of $R$. We will also discuss a version of the above theorem for the case when the residue field is finite at the end of \S~\ref{sectionmain}. The above theorem allows us to find sufficiently many local normal domains of mixed characteristic as specializations, but it does not tell us how to find $U$. We will show how to find $U$ in Example \ref{algorithm}. It is worth pointing out that if $\dim R=2$, the local Bertini theorem fails due to a simple reason. In fact, if the quotient $R/yR$ is normal, then it is a discrete valuation ring, so $R$ must be regular. By Cohen's structure theorem, there is a surjection $A[[z_0,\ldots,z_d]] \twoheadrightarrow R$, where $d+1$ is the number of the minimal generators of $\fm$, so that the minimal generators of $\fm$ are just the image of $z_0,z_1,\ldots,z_d$ under this surjection. Related to our main result, we mention that if $(R,\fm)$ is a local ring and $y \in \fm$ is a non-zero divisor such that $R/yR$ is normal, then $R$ is normal (\cite{EGA} 5.12.7).

In \S~\ref{sectionsubmain}, we establish a fact that the parameter set of specializations in the local Bertini theorem is infinite (see Proposition \ref{prop2}), which is used to prove the main result on characteristic ideals.

In \S~\ref{sectionmainSerre}, we prove a version of the local Bertini theorem in the case when Serre's condition $(R_n)$ or $(S_n)$ is satisfied.

In \S~\ref{characteristic}, we define characteristic ideals of finitely generated torsion modules over a Noetherian normal domain as reflexive ideals, following the paper \cite{SkUr}. 

In \S~\ref{application}, the main purpose is to generalize the results proved in \cite{Oc1} over a regular ring to the case over a Noetherian normal domain. First, we prove some preliminary results on characteristic ideals. Then combining Main Theorem A, we prove another main theorem (see Theorem \ref{charid}) as follows:

\begin{maintheoremb}[Control Theorem for Characteristic Ideals]
Let $(R,\fm,\mathbf{k})$ be a complete local domain of mixed characteristic $p>0$ satisfying the conditions (1), (2) and (3) of Main Theorem A, except that we now assume the residue field to be finite. Let $M$ and $N$ be finitely generated torsion $R$-modules. Then for certain infinite subsets $\mathcal{L}_{W(\overline{\mathbb{F}})}(M_{W(\overline{\mathbb{F}})})$ and $\mathcal{L}_{W(\overline{\mathbb{F}})}(N_{W(\overline{\mathbb{F}})})$ of the 
projective space $\mathbb{P}^{d}(W(\overline{\mathbb{F}}))$ defined in a natural way (cf. Definition \ref{definition:L}), the following statements are equivalent:

\begin{enumerate}
\item[$\mathrm{(1)}$]
$\Char_R(M) \subseteq \Char_R(N)$, where $\Char_R(-)$ denotes the characteristic ideal (cf. Definition \ref{definition:characteristicideal}).

\item[$\mathrm{(2)}$]
For all but finitely many height-one primes:
$$
\mathbf{x}_{\widetilde{a}}R_{W(\overline{\mathbb{F}})} \in  \mathcal{L}_{W(\overline{\mathbb{F}})}(M_{W(\overline{\mathbb{F}})}) \cap \mathcal{L}_{W(\overline{\mathbb{F}})}(N_{W(\overline{\mathbb{F}})}),
$$
we have
$$
\Char_{R_{W(\mathbb{F}')}/\mathbf{x}_{\widetilde{a}}R_{W(\mathbb{F}')}}(M_{W(\mathbb{F}')}/\mathbf{x}_{\widetilde{a}}M_{W(\mathbb{F}')}) \subseteq \Char_{R_{W(\mathbb{F}')}/\mathbf{x}_{\widetilde{a}}R_{W(\mathbb{F}')}}(N_{W(\mathbb{F}')}/\mathbf{x}_{\widetilde{a}}N_{W(\mathbb{F}')}),
$$
where $\mathbb{F}'$ is any finite field extension of $\mathbb{F}$ depending on ${\widetilde{a}}$ such that $\mathbf{x}_{\widetilde{a}} \in R_{W(\mathbb{F}')}$.

\item[$\mathrm{(3)}$]
For all but finitely many height-one primes:
$$
\mathbf{x}_{\widetilde{a}}R_{W(\overline{\mathbb{F}})} \in  \mathcal{L}_{W(\overline{\mathbb{F}})}(M_{W(\overline{\mathbb{F}})}) \cap \mathcal{L}_{W(\overline{\mathbb{F}})}(N_{W(\overline{\mathbb{F}})}),
$$
we have
$$
\Char_{R_{W(\overline{\mathbb{F}})}/\mathbf{x}_{\widetilde{a}}R_{W(\overline{\mathbb{F}})}}(M_{W(\overline{\mathbb{F}})}/\mathbf{x}_{\widetilde{a}}M_{W(\overline{\mathbb{F}})}) \subseteq \Char_{R_{W(\overline{\mathbb{F}})}/\mathbf{x}_{\widetilde{a}}R_{W(\overline{\mathbb{F}})}}(N_{W(\overline{\mathbb{F}})}/\mathbf{x}_{\widetilde{a}}N_{W(\overline{\mathbb{F}})}).
$$
\end{enumerate} 
\end{maintheoremb}

Main Theorem B will be crucial in a forthcoming paper \cite{OcSh}, where we plan to compare characteristic ideals of certain torsion modules arising from Iwasawa theory as developed in \cite{Oc1} (for example, those torsion modules arise as the Pontryagin dual modules of Selmer groups associated with two-dimensional certain Galois representations with values in a complete local ring with finite residue field).

\section{Specialization map and the topology for Bertini-type theorems}
\label{topology}

In this section, we introduce notation and make definitions. We fix a prime integer $p>0$. First, we discuss the specialization map. Let $(A,\pi_A,\mathbf{k})$ be a complete discrete valuation ring and let us choose a point:
$$
a=(a_0:\cdots:a_n) \in \mathbb{P}^n(A).
$$ 
Suppose that the valuation $v$ of $\pi_A$ is equal to 1 and that the $i$-th coordinate $a_i$ has the minimal valuation $v_i:=v(a_i)$ among the valuations of $a_0, \ldots ,a_n$. We put $a'_j =a_j /\pi^{v_i}_A$ 
($0\leq j \leq n$). Then we have a presentation $a=(a'_0:\cdots:a'_n)$, so that $a_j \in A$ ($0\leq j\leq n$) and that $a'_i$ is a unit of $A$. So we may and we will 
think of the projective space as 
$$
\mathbb{P}^n(A)=\{\text{homogeneous coordinates }(a_0:\cdots:a_n)\text{ in $A$} ~\vert~\mbox{$a_i \in A^\times$ for some $i$}\} / \sim 
$$
where the equivalence relation $\sim$ is a simultaneous multiplication by a unit of $A$.

\begin{definition}
\label{definition:basic_p_adic}
Let $(A,\pi_A,\mathbf{k})$ be a complete discrete valuation ring and pick a point $a = (a_0 :\cdots:a_n ) \in \mathbb{P}^n(A)$.
\begin{enumerate}
\item
We define a \textit{lift} $\widetilde{a}= (\widetilde{a}_0,\ldots,\widetilde{a}_n ) \in \mathbb{A}^{n+1}(A)$ of $a \in \mathbb{P}^n(A)$ to be an element belonging to the set consisting of all inverse images of $a \in \mathbb{P}^n(A)$ via the projection map: 
$$
\mathbb{A}^{n+1} (A) \setminus \{\text{points whose coordinates are all non-units}\} \twoheadrightarrow \mathbb{P}^n(A).
$$
 
\item 
Let $\widetilde{a} \in \mathbb{A}^{n+1}(A)$ be a lift of $a \in \mathbb{P}^n(A)$. Then $\widetilde{a}$ gets mapped to a point $\overline{a}=(\overline{a}_0:\cdots:\overline{a}_n) \in \mathbb{P}^n(\mathbf{k})$ via the surjection map $A \to \mathbf{k}:=A/\pi_{A} A$. Thus, we construct a \textit{specialization map}:
$$ 
\Sp_A:\mathbb{P}^n(A) \to \mathbb{P}^n(\mathbf{k})
$$
by setting $\Sp_A(a)=\overline{a}$. This map does not depend on the choice of the lift of $a$.
\end{enumerate}
\end{definition}

The set $\mathbb{P}^n(\mathbf{k})$ is endowed with the Zariski topology, while $\mathbb{P}^n(A)$ is endowed with the topology induced by the valuation on $A$. Hence we simply regard $\mathbb{P}^n(A)$ as a set of points equipped with this topology. It is also straightforward from the definition that we naturally identify the set $\mathbb{P}^n (A)$ with $\mathbb{P}^n (\mathrm{Frac} (A))$ 
as well as the induced topology.  

We begin to pin down the suitable topology for formulating Bertini-type theorems in mixed characteristic. Let $(R,\fm,\mathbf{k})$ be a local Noetherian ring. Then we say that a reduced local ring $(R,\fm,\mathbf{k})$ is of \textit{mixed characteristic $p>0$}, if every component of the total ring of fractions of $R$ is of characteristic zero and the residue field $\mathbf{k}$ is of characteristic $p$. 

Now we assume that $(R,\fm,\mathbf{k})$ is a complete reduced local ring of mixed characteristic $p>0$ and $(A,\pi_A,\mathbf{k})$ is a complete discrete valuation ring such that $\pi_A A=pA$ and there is an injection $A \hookrightarrow R$ of rings, which induces an isomorphism on residue fields, say $\mathbf{k}=A/\pi_A A \simeq R/\fm$. In this situation, we call $(A,\pi_A,\mathbf{k})$ together with a map $A \hookrightarrow R$ a \textit{coefficient ring} of $R$. 

\begin{example}
Let $R:=\mathbb{Z}_p[[x,y]]/(p-xy)$. Then $R$ is a finite extension of $\mathbb{Z}_p[[x+y]]$ defined by the Eisenstein equation $t^2-(x+y)t+p=0$ and $\mathbb{Z}_p$ is a coefficient ring of $R$.
\end{example}

In what follows, we will fix a coefficient ring $(A,\pi_A,\mathbf{k})$. We denote by $\mathbf{Loc.alg}_{/A}$ the category of complete local $A$-algebras for a discrete valuation ring $A$. Note that the category of complete local $\mathbf{k}$-algebras is a full subcategory of $\mathbf{Loc.alg}_{/A}$. Let $(R,\fm,\mathbf{k})$ be a local Noetherian ring and let $M$ be a finitely generated $R$-module. We say that a set of elements $x_0,\ldots,x_d$ of $M$ is the \textit{minimal generators} of $M$, if the following conditions hold: 
$$
\sum_{i=0}^d \mathbf{k} \cdot \overline{x}_i=M/\fm M~\mbox{and}~\dim_{\mathbf{k}} M/\fm M=d+1,
$$ 
where $\overline{x}_i$ denotes the image of $x_i$ in $M/\fm M$.

\begin{definition}
\label{definition:Bertini_property}
Let $(R,\fm,\mathbf{k}) \in \mathbf{Loc.alg}_{/A}$ be a reduced local ring of mixed characteristic $p>0$.

\begin{enumerate}
\item 
Fix a set of minimal generators $x_0 ,\ldots ,x_n \in \fm$ and let $a =(a_0:\cdots:a_n) \in \mathbb{P}^n(A)$. For any lift $\widetilde{a} = (\widetilde{a}_0 ,\ldots,\widetilde{a}_n ) \in \mathbb{A}^{n+1}(A)$ of $a$ in the sense of Definition \ref{definition:basic_p_adic}, we put
$$
\mathbf{x}_{\widetilde{a}}:=\sum_{i=0}^{d} \widetilde{a}_{i}x_{i}. 
$$ 
By definition, the element $\mathbf{x}_{\widetilde{a}} \in R$ depends on the choice of a lift $\widetilde{a}$, but the ideal $\mathbf{x}_{\widetilde{a}} R$ depends only on $a \in \mathbb{P}^n(A)$. 

\item
With the notation as above, assume that $\mathbf{P}$ is a ring-theoretic property on Noetherian rings. 
We say that $(R,\fm,\mathbf{k}) \in \mathbf{Loc.alg}_{/A}$ satisfies a \textit{local Bertini theorem} for the property $\mathbf{P}$ if for a fixed set of minimal generators $x_0 ,\ldots ,x_n \in \fm$, 
there exists a Zariski (dense) open subset $U \subseteq \mathbb{P}^n(\mathbf{k})$ such that $R/\mathbf{x}_{\widetilde{a}}R$ has $\mathbf{P}$ for all $a=(a_0:\cdots :a_n) \in \Sp_A^{-1}(U) \subseteq \mathbb{P}^n(A)$.
\end{enumerate}
\end{definition}

One may formulate the local Bertini theorem in a different way. For example, the completeness of $R$ may be dispensed with, or a set of minimal generators of the maximal ideal of $R$ may be replaced with a set of minimal generators of a smaller ideal. In this article, we adopt the above definition.

\begin{remark}
The naturality of the above definition is explained as follows. We endow $\mathbb{P}^n(\mathbf{k})$ with the Zariski topology. Let $\overline{f} \in \mathbf{k}[x_0,\ldots,x_n]$ be a non-zero homogeneous polynomial, let $f \in A[x_0,\ldots,x_n]$ be any fixed homogeneous lifting of $\overline{f}$, and let $U^+(\overline{f}) \subseteq \mathbb{P}^n(\mathbf{k})$ be an open subset defined by $\overline{f} \ne 0$. Then the inverse image of the open subset $U^+(\overline{f})$ under the map $\Sp_A:\mathbb{P}^n(A) \to \mathbb{P}^n(\mathbf{k})$ can be described as follows:
$$
\Sp_A^{-1}(U^+(\overline{f}))=\{a=(a_0:\cdots:a_n) \in \mathbb{P}^n(A)~|~f(a) \in A^{\times}\}.
$$
This is an admissible open subset of $\mathbb{P}^n(A)$ (\cite{FrPut} for this fact). Our objective is to show that this topology is suitable in formulating the local Bertini theorem in the mixed characteristic case.
\end{remark}

The following proposition is indispensable for the proof of Theorem \ref{theorem2} and Theorem \ref{thm:localBertini}.

\begin{proposition}
\label{prop1}
Let $U \subseteq \mathbb{A}^n(L)$ be any non-empty Zariski open subset for an infinite field $L$. Then $U$ is dense. Furthermore, if $K \subseteq L$ is any subfield of $L$ such that $K$ is infinite, the intersection $U \cap \mathbb{A}^n(K)$ is also a Zariski dense open subset of $\mathbb{A}^n(K)$. The above assertions hold over the projective space as well.
\end{proposition}

\begin{proof}
The first assertion about density is obvious. So we prove the second assertion. Then it suffices to prove the following statement:
\begin{enumerate}
\item[$\bullet$]
Let $V:=V(I)$ be a Zariski closed subset of $\mathbb{A}^n_L:=\Spec L[X_1,\ldots,X_n]$, where $I$ is an ideal of $L[X_1,\ldots,X_n]$.
If $\mathbb{A}^n(K)$ is contained in the set of $K$-rational points of $V$, then $V$ is equal to $\mathbb{A}^n_L$.
\end{enumerate}
Let us prove this. If $\mathbb{A}^n(K)$ is contained in the set of $K$-rational points of $V$, we have
$$
I \subseteq \bigcap_{(a_1,\ldots,a_n) \in K^n} (X_1-a_1,X_2-a_2,\ldots,X_n-a_n) \subseteq L[X_1,\ldots,X_n].
$$
Since $K$ is infinite, we have
$$
\bigcap_{(a_1,\ldots,a_n) \in K^n} (X_1-a_1,X_2-a_2,\ldots,X_n-a_n)=0
$$
and thus $I=(0)$. This implies that $V(I)=\mathbb{A}^n_L$.
\end{proof}

\section{Discussion on basic elements}
\label{basicelement}

We start with the definition of basic elements.

\begin{definition}
Let $R$ be a Noetherian ring and let $M$ be a finitely generated $R$-module. 
\begin{list}{}{}
\item[(1)] 
Assume that $(R,\fm)$ is local. The depth of $M$, denoted by $\depth_{R}M$, is defined to be the maximal length of all $M$-regular sequences contained in $\fm$. 
\item[(2)]
Let $\fp \in \Spec R$. Then $\mu_{\fp}(M)$ denotes the number of minimal generators of the $R_{\fp}$-module $M_{\fp}$.  
\item[(3) (Swan)]
Let $\fp \in \Spec R$. We say that an element $m \in M$ is \textit{basic} at $\fp$, if we have
$$
\mu_{\fp}(M)-\mu_{\fp}(M/(R \cdot m))=1.
$$
More generally, we say that a set of elements $m_1,\ldots,m_n$ of $M$ is \textit{k-fold basic} at $\fp$, if the following inequality holds:
$$
\mu_{\fp}(M)-\mu_{\fp}(M/(\sum_{i=1}^{n} R \cdot m_{i})) \ge k,
$$
that is, $N=R \cdot m_1+\cdots+R\cdot m_n$ contains at least $k$ minimal generators at $\fp$.
\end{list}
\end{definition}

\begin{remark}\label{remark:32}
Under the setting of the above definition, 
let $M^{(r)}:=M/(\sum_{i=1}^{r} R \cdot m_{i})$ for a set of elements $m_1,\ldots,m_k$ of $M$ and $r$ satisfying $0 \le r \le k-1$ and pick a prime ideal $\fp$ of $R$. Then  
$$
\mu_{\fp}(M^{(r)})-\mu_{\fp}(M^{(r)}/(R\cdot m_{r+1}))=1 \iff m_{r+1} \notin \fp M^{(r)}_{\fp}
$$
for $0 \le r \le k-1$ by Nakayama's lemma. In other words, $m_1,\ldots,m_k$ form partial generators of the $k(\fp)$-vector space $M \otimes_R k(\fp)$ with $k(\fp):=R_{\fp}/\fp R_{\fp}$.
\end{remark}

We shall use (finite) K\"ahler differentials (\cite{Kunz} as a reference). For a complete local ring $(R,\fm)$ with its coefficient ring $A$, the usual module of K\"ahler differentials $\Omega_{R/A}$ is not a finitely generated $R$-module. Instead, one uses the completed module $\widehat{\Omega}_{R/A}$. This is the $\fm$-adic completion of $\Omega_{R/A}$ and it is a finitely generated $R$-module. It can be also defined as follows. Let $I$ denote the kernel of the map $\mu:R \widehat{\otimes}_{A} R \to R$ defined by $\mu(a \otimes b)=ab$. Then $\widehat{\Omega}_{R/A}:=I/I^2$. Let $d: R \to \widehat{\Omega}_{R/A}$ be the canonical derivation defined by $a \mapsto a \otimes 1-1 \otimes a$. The connection of K\"ahler differential modules with the symbolic power ideals is expressed by the following simple fact (\cite{Fle}, Lemma 2.2 for its proof).

\begin{lemma}
\label{lemma1}
Let $M$ be a module over a ring $R$, let $\fp$ be a prime of $R$, and let $d:R \to M$ be a derivation. If for $x \in R$, $dx \in M$ is basic at $\fp$, then $x \notin \fp^{(2)}$.
\end{lemma}

Let $I \subseteq R$ be an ideal. We denote by $\Min_R(I)$ the set of all prime ideals that are minimal over $I$. The authors are grateful to Prof. V. Trivedi for explaining the proof of the following lemma (\cite{Fle}, Lemma 1.2).

\begin{lemma}[Flenner]
\label{lemma2}
Suppose that $R$ is a Noetherian ring, $M$ is a finitely generated $R$-module, $U \subseteq \Spec R$ is a Zariski open subset, and $\{m_0,\ldots,m_n\}$ is a set of elements of $M$, which generates the submodule $N \subseteq M$. Suppose that we have $t\in \mathbb{Z}$ $($which can be negative$)$ such that  
$$
\mu_{\fp}(M)-\mu_{\fp}(M/N) \ge \dim(V(\fp ) \cap U) -t
$$
for every $\fp \in U$. Let $(\phi^*)^{-1}(U)$ be the inverse image of $U$ under $\phi^*:\Spec R[X_0,\ldots,X_n] \to \Spec R$ induced by the natural injection $\phi:R \to R[X_0,\ldots,X_n]$. Then there exists an ideal $(F_1,\ldots,F_r) \subseteq R[X_0,\ldots,X_n]$ such that
$$
\dim (V(F_1,\ldots,F_r) \cap (\phi^*)^{-1}(U)) \le n+1+t
$$ 
and the element
$$
\sum_{i=0}^{n} m_{i} \otimes X_{i} \in M \otimes_{R} R[X_0,\ldots,X_n]
$$
is basic on $U(F_1,\ldots,F_r) \cap (\phi^*)^{-1}(U)$.
\end{lemma}

We need the following technical lemma for the proof of the main theorem. We recall that a local domain $S$ is \textit{catenary}, if and only if $\Ht \fp+\dim S/\fp=\dim S$ for all $\fp \in \Spec S$.

\begin{lemma}
\label{prelim}
Let $(R,\fm)$ be an excellent local domain and let $\phi:R \to R[X_0,\ldots,X_d]$ be a natural injection with an ideal $(F_1,\ldots,F_r) \subseteq R[X_0,\ldots,X_d]$. Suppose that the following conditions hold:

\begin{enumerate}
\item[$\mathrm{(1)}$]
$\dim \big(V(F_1,\ldots,F_r) \cap (\phi^*)^{-1}(U(\fm))\big) \le d+1$.

\item[$\mathrm{(2)}$]
$(F_1,\ldots,F_r) \subseteq \fm R[X_0,\ldots,X_d]$.
\end{enumerate}
Then the set of prime ideals $P \in \Spec R[X_0,\ldots,X_d]$ satisfying the property:
$$
(F_1,\ldots,F_r) \subseteq P \subseteq \fm R[X_0,\ldots,X_d]
$$ 
is finite.
\end{lemma}

\begin{proof}
Let $I:=(F_1,\ldots,F_r)$ and consider the finite set:
$$
S:=\{P \in \Min_{R[X_0,\ldots,X_d]}(I)~|~P \subseteq \fm R[X_0,\ldots,X_d]\}.
$$
We would like to show that
\begin{equation}
\label{chain-1}
0 \le \Ht(\fm R[X_0,\ldots,X_d])-\Ht P \le 1
\end{equation}
for $P \in S$. If $(\ref{chain-1})$ is true for $P \in S$, then it is immediate that a prime ideal between $P$ and $\fm R[X_0,\ldots,X_d]$ is either $P$ or $\fm R[X_0,\ldots,X_d]$. Hence if $(\ref{chain-1})$ is proved for all $P \in S$, the lemma follows. So we will prove $(\ref{chain-1})$ for all $P \in S$.

Since $R$ is an excellent local domain, $R[X_0,\ldots,X_d]$ is a catenary domain. From this it follows that 
$0 \le \Ht(\fm R[X_0,\ldots,X_d])-\Ht P$ for $P \in S$. If we have $\Ht(\fm R[X_0,\ldots,X_d])=\Ht P$, then $\fm R[X_0,\ldots,X_d]=P$. So it suffices to establish the inequality:
$$
\Ht(\fm R[X_0,\ldots,X_d])-\Ht P \le 1.
$$
Let $T$ be the localization of $R[X_0,\ldots,X_d]$ at the prime ideal $\fm R[X_0,\ldots,X_d]$. Then $R \to T$ is a flat local map of local domains and $\dim R=\dim T$. Now assume that $P \ne \fm R[X_0,\ldots,X_d]$ for $P \in S$. Since $R$ is an excellent local domain, $T$ is a catenary local domain and we have $\Ht P=\dim T-\dim T/PT$. So if we can prove
$$
\dim T/PT=1,
$$
the relation $(\ref{chain-1})$ follows from this and we are done. Let 
\begin{equation}
\label{chain0}
P_0:=P \subsetneq P_1 \subsetneq \cdots \subsetneq P_s=\fm R[X_0,\ldots,X_d]
\end{equation}
 be a chain of prime ideals of maximal length. Since $\dim T/PT=s$, it suffices to prove $s=1$. Consider a chain of prime ideals ($0 \le k \le s$) in $R[X_0,\ldots,X_d]$:
\begin{equation}
\label{chain1}
P_0 \subsetneq \cdots \subsetneq P_k \subsetneq P_k+(X_0) \subsetneq P_k+(X_0,X_1) \subsetneq \cdots \subsetneq P_k+(X_0,\ldots,X_d).
\end{equation}
On the other hand, $\dim (V(F_1,\ldots,F_r) \cap (\phi^*)^{-1}(U(\fm)))$ is equal to the length of the chain of prime ideals of maximal length in $R[X_0,\ldots,X_d]$:
\begin{equation}
\label{chain2}
Q_0 \subsetneq Q_1 \subsetneq \cdots \subsetneq Q_t~\mbox{such that}~I \subseteq Q_0~\mbox{and}~\phi^{-1}(Q_t) \ne \fm.
\end{equation}

First, note that $s=0$ is impossible, since we assumed 
$P \ne \fm R[X_0,\ldots,X_d]$. So in the rest of the proof, let us assume $s \ge 2$ and derive 
a contradiction. In this case, taking $k=1$ in $(\ref{chain1})$, we get a chain of prime ideals of length $d+2$:
\begin{equation}
\label{chain3}
P_0 \subsetneq P_1 \subsetneq P_1+(X_0) \subsetneq P_1+(X_0,X_1) \subsetneq \cdots \subsetneq P_1+(X_0,\ldots,X_d).
\end{equation}
Since $P_1 \subsetneq P_2 \subseteq \fm R[X_0,\ldots,X_d]$ by $s \ge 2$, we have $\phi^{-1}(P_1+(X_0,\ldots,X_d)) \ne \fm$. By comparing $(\ref{chain2})$ and $(\ref{chain3})$, we must have $d+2 \le t$, because $(\ref{chain2})$ is of maximal length. However, this is not compatible with the condition $(1)$. Therefore, we get $s=1$ as desired.
\end{proof}

In this paper, the Noetherian induction and the lemma of generic freeness are important tools. For this, we often need the following fact and use it freely.

\begin{lemma}
\label{lemma:noetherian_induction}
Let $R$ be a Noetherian ring and let $M$ be a finitely generated $R$-module. 
Let $I$ be an ideal of $R$ and let $\mathfrak{p}$ be a prime ideal of $R$ satisfying $I \subseteq \mathfrak{p}$. We denote the quotient ring $R /I$ by $\overline{R}$ 
and the image of $\mathfrak{p}$ in $\overline{R}$ by $\overline{\mathfrak{p}}$. 
 Then an element $x \in M$ is basic at $\fp \in \Spec R$ if and only if $\overline{x} \in \overline{M}:=M/IM$ is basic at $\overline{\fp} \in \Spec \overline{R}$.
\end{lemma}

\begin{proof} 
Note that there is a commutative square:
$$
\begin{CD}
M @>>> \overline{M} \\
@V\pi VV @V\overline{\pi}VV \\
M \otimes_R k(\fp) @=\overline{M} \otimes_{\overline{R}} k(\overline{\fp})
\end{CD}
$$
In this diagram, we have $\pi(x)=\overline{\pi}(\overline{x})$. By Remark \ref{remark:32}, $x \in M$ (resp. $\overline{x} \in \overline{M}$) is basic at $\fp \in \Spec R$ (resp. $\overline{\fp} \in \Spec \overline{R}$) if and only if $\pi(x) \in M \otimes_R k(\fp)$ (resp. $\overline{\pi}(\overline{x}) \in \overline{M} \otimes_{\overline{R}} k(\overline{\fp})$) is not zero. This completes the proof.
\end{proof}

Let $R$ be a Noetherian ring and let $M$ be a finitely generated $R$-module. Fix an element $x \in M$. Let us define
$$
Z_x:=\{\fp \in \Spec R~|~x~\mbox{is basic at}~\fp\}.
$$
Then $Z_x$ is not necessarily a Zariski open subset. However, we have the following result.

\begin{proposition}
Let $R$ be a Noetherian ring and suppose $M$ is a finitely generated $R$-module. Fix an element $x \in M$ and let $Z_x$ be as above. Then $Z_x$ is a $($possibly empty$)$ constructible subset of $\Spec R$.
\end{proposition}

\begin{proof}
First we consider the special case where $M$ is finitely generated and projective over $R$. 
Recall that for a finitely generated projective $R$-module $N$, $x \in N$ is basic at $\fp \in \Spec R$ if and only if $R_{\fp} \cdot x$ spans a direct summand of 
$N_{\fp}$ (see \cite{EisEv}, Lemma 1 for the proof). 
Note that, for any short exact sequence of $R$-modules, the locus on which it splits is a Zariski open subset of $\mathrm{Spec} R$. Hence, we proved that $Z_x$ is Zariski open when $M$ is finitely generated and projective over $R$. 

Then we will prove the case for general $R$-modules $M$. First, we apply 
Lemma \ref{lemma:noetherian_induction} by taking $I$ to be 
the nil-radical $\sqrt{0}$ of $R$ and by taking $\overline{R}$ 
to be $R_{\red}:=R/\sqrt{0}$. Since $R \twoheadrightarrow R_{\red}$ 
induces an isomorphism of topological spaces 
$\Spec R_{\red} \overset{\sim}{\longrightarrow} \Spec R$, 
we may and we will assume that $R$ is reduced from now on.  
Thus there exists $f \in R$ for which $M[f^{-1}]$ is a free $R[f^{-1}]$-module. Denote by $U \subseteq \Spec R[f^{-1}]$ the locus on which $x \in M[f^{-1}]$ is basic. 

Continuing this process, we may find a chain of closed subsets:
$$
V(f)=:V_1 \supseteq V_2 \supseteq V_3 \supseteq \cdots 
$$
such that $x \in M$ is basic at $\fp \in V_i$ if and only if $\fp \in Z_i:=V_i\backslash V_{i-1} \subseteq V_i$. Moreover, $Z_i$ is an open subset of $V_i$. Since $\Spec R$ is a Noetherian space, this chain stabilizes. Thus, there is an integer $N>0$ such that $V_N=V_{N+1}=\cdots$. We set $Z:=Z_1 \cup Z_2 \cup \cdots \cup Z_N$. Then it follows that $x \in M$ is basic at each point of $Z$ and $Z$ is a constructible subset of $\Spec R$. So $Z_x:=U \cup Z$ is the sought one and constructible.
\end{proof}

\section{Main Theorems}\label{sectionmain}

In this section, we establish our main theorems.

\begin{lemma}
\label{lemma3}
Let $(A,\pi_A,\mathbf{k})$ be a discrete valuation ring and let $f \in A[y_1,\ldots,y_d]$ be a non-zero $($possibly constant$)$ polynomial. Then there exists $t \in \mathbb{Z}_{\ge 0}$ such that $\pi_A^{-t}f \in A[y_1,\ldots,y_d]$ and the reduction of $\pi_A^{-t}f$ modulo $\pi_A$ is a non-zero $($possibly constant$)$ polynomial in $\mathbf{k}[y_1,\ldots,y_d]$.
\end{lemma}

\begin{proof}
The proof goes by induction on $d$. If $d=1$, we write $f=a_my^m_1+a_{m-1}y^{m-1}_1+\cdots+a_0$ with $a_i \in A$. Let $0 \le h \le m$ be such that the valuation $v(a_h)$ is the smallest in the set $\{v(a_i)~|~0 \le i \le m\}$ and write $a_h=(\mbox{unit})\cdot \pi_A^t$. Dividing $f$ by $\pi_A^t$, we get the desired polynomial.

In general, write $f=b_ny^n_d+b_{n-1}y^{n-1}_d+\cdots+b_0$ for $b_i \in A[y_{1},\ldots,y_{d-1}]$. Applying the induction hypothesis to every $b_i$, we may find $t_i \in \mathbb{Z}_{\ge 0}$ such that $\pi_A^{-t_i} b_i$ has the desired property. Let $t_s:=\min\{t_i~|~0 \le i \le n\}$. Then the term $\pi_A^{-t_s}b_sy^s_d$ modulo $\pi_A$ is non-zero and it is clear that $\pi_A^{-t_s}f$ is contained in $A[y_1,\ldots,y_d]$. Hence it suffices to put $t:=t_s$.
\end{proof}

The following lemma will play a role in the final step of the proof of the next theorem.

\begin{lemma}
\label{Case1Case2}
Let $(R,\fm,\mathbf{k})$ be a complete local domain of mixed characteristic $p > 0$ with infinite residue field $\mathbf{k}$ and a coefficient ring $(A,\pi_A)$. Fix a set of minimal generators $x_0,x_1,\ldots,x_d$ of $\fm$ together with a prime ideal $\fp$ of $R$ with $\fp \ne \fm$. Then there exists a non-empty open subset $U \subseteq \mathbb{P}^d(\mathbf{k})$ such that
$$
\mathbf{x}_{\widetilde{a}}:=\sum_{i=0}^d \widetilde{a}_ix_i \notin \fp
$$
for every $a=(a_0:\cdots:a_d) \in \Sp_A^{-1}(U)$.
\end{lemma}

\begin{proof}
Consider the homogeneous polynomial:
$$
F(X_0,\ldots,X_d):=\sum_{i=0}^dx_iX_i \in R[X_0,\ldots,X_d].
$$
Then we have an equivalence of conditions:
$$
\mathbf{x}_{\widetilde{a}}=\sum_{i=0}^{d} \widetilde{a}_i x_i \notin \fp \iff F(\widetilde{a}_{0},\ldots,\widetilde{a}_{d}) \not\equiv 0 \pmod {\fp}.
$$
We will divide the proof of the lemma into the separate cases.

${\bf{Case 1}}$: 
Assume that $\pi_A \in \fp$. Let $S$ be the localization of $R$ at $\fp$. Then $A \to S$ is a flat local map of local rings and $S$ is of mixed characteristic. Let $\mathbf{k}'$ denote the residue field of $S$. Then we have a mapping:
$$
\mathbb{P}^d(\mathbf{k}) \to \mathbb{P}^d(\mathbf{k}').
$$ 
Since $F$ is a homogeneous polynomial, the condition $F(\widetilde{a}_{0},\ldots,\widetilde{a}_{d}) \not\equiv 0 \pmod {\fp}$ defines an open subset $V \subseteq \mathbb{P}^d(\mathbf{k}')$. By Proposition \ref{prop1}, 
$$
U:=V \cap \mathbb{P}^d(\mathbf{k})
$$ 
is a dense open subset of $\mathbb{P}^d(\mathbf{k})$, since the field $\mathbf{k}$ is infinite. Then this $U$ is a non-empty open subset of $\mathbb{P}^d(\mathbf{k})$ with the desired property.

${\bf{Case 2}}$:
Assume that $\pi_A \notin \fp$. In this case, notation being as in ${\bf{Case 1}}$, we have a mapping:
$$
\mathbb{P}^d (\mathrm{Frac}(A)) \to \mathbb{P}^d(\mathbf{k}').
$$
Again, the condition that $F(\widetilde{a}_{0},\ldots,\widetilde{a}_{d}) \not\equiv 0 \pmod {\fp}$ defines a Zariski open subset $V \subseteq \mathbb{P}^d(\mathbf{k}')$. We use an identification $\mathbb{P}^d (\mathrm{Frac}(A))=\mathbb{P}^d (A)$.

Proposition \ref{prop1} implies that $V \cap \mathbb{P}^d( \mathrm{Frac}(A))$ is a dense open subset of $\mathbb{P}^d( \mathrm{Frac}(A))$, and it is covered by open sets of the form $U^+(f)$
for a homogeneous polynomial $f \in A[X_0,\ldots,X_d]$. Let $\overline{f} \in \mathbf{k}[X_0,\ldots,X_d]$ be its reduction modulo $\pi_A$ and let $U^+(\overline{f}) \subseteq \mathbb{P}^d(\mathbf{k})$ be the corresponding open subset. In view of Lemma \ref{lemma3}, by replacing $f$ with $\pi_A^{-t}f$ if necessary, we may assume that $f \in A[X_0,\ldots,X_d]$ and $\overline{f} \in \mathbf{k}[X_0,\ldots,X_d]$ is not equal to 0. Since $f(\widetilde{a}) \in A^{\times}$ if and only if $\overline{f}(\Sp_A(a)) \ne 0$,
we have
$$
\Sp_A^{-1}(U^+(\overline{f}))=\{a=(a_0:\cdots:a_d) \in \mathbb{P}^d(A)~|~f(\widetilde{a}) \in A^{\times}\} \subseteq \mathbb{P}^d(A).
$$ 
Applying this description to the open covering of $V \cap \mathbb{P}^d(\mathrm{Frac}(A))$ consisting of $U^+(f)$, the image 
$$
U:=\Sp_A(V \cap \mathbb{P}^d(\mathrm{Frac}(A))) \subseteq \mathbb{P}^d(\mathbf{k})
$$ 
is a Zariski open subset. Then this $U$ is a non-empty open subset of $\mathbb{P}^d(\mathbf{k})$ with the desired property.

Combining ${\bf{Case 1}}$ and ${\bf{Case 2}}$ together, we finish the proof of the lemma.
\end{proof}

The following is our first main theorem.

\begin{theorem}
\label{theorem2}
Let $(R,\fm,\mathbf{k})$ be a complete local domain of mixed characteristic $p > 0$ and suppose the following conditions:

\begin{enumerate}
\item[$\mathrm{(1)}$]
let $A \to R$ be a coefficient ring map for a complete discrete valuation ring $(A,\pi_{A})$; 

\item[$\mathrm{(2)}$]
let $x_0,x_1,\ldots,x_d$ be a fixed set of minimal generators of $\fm$;

\item[$\mathrm{(3)}$]
the residue field $\mathbf{k}$ is infinite.
\end{enumerate}
Then there exists a Zariski dense open subset $U' \subseteq \mathbb{P}^d(\mathbf{k})$ such that we have 
$$
\mathbf{x}_{\widetilde{a}}=\sum_{i=0}^{d} \widetilde{a}_i x_{i} \notin \fp^{(2)}
$$
for every prime $\fp$ of $R$ and every $a=(a_0:\cdots:a_d) \in \Sp^{-1}_A(U') \subseteq \mathbb{P}^d(A)$. 
\end{theorem}

Note that $x_0,x_1,\ldots,x_d$ are not a system of parameters of $R$, unless $R$ is a regular local ring.

\begin{proof}
First, the $R$-module $\widehat{\Omega}_{R/A}$ is generated by $dx_0,\ldots,dx_d$, which is easily verified by considering the surjective ring map $A[[X_0,\ldots,X_d]] \twoheadrightarrow R$ defined by mapping each $X_i$ to $x_i$. Let us recall an important fact. If $R$ is a complete local domain with $(A,\pi_A)$ its coefficient ring, it follows from (\cite{Fle}, Lemma 2.6, or \cite{Triv}, Lemma 2) that for any fixed element $\fp \in U(\fm)$, we have
$$
\mu_{\fp}(\widehat{\Omega}_{R/A}) \ge \dim(R/\fp)-1.
$$
Noting that $\dim(R/\fp)-1=\dim (V(\fp) \cap U(\fm))$ for a subscheme $V(\fp) \cap U(\fm) \subseteq \Spec R$ and applying Lemma \ref{lemma2} for the $R$-module $\widehat{\Omega}_{R/A}$, there is an ideal $(F_1,\ldots,F_r) \subseteq R[X_{0},\ldots,X_{d}]$ such that we have 
\begin{equation}
\label{zero}
\dim (V(F_1,\ldots,F_r ) \cap (\phi^*)^{-1}(U(\fm)) ) \le d+1 ,
\end{equation}
where $(\phi^*)^{-1}(U(\fm))$ is the inverse image of $U(\fm)$ under $\phi^*:\Spec R[X_0,\ldots,X_d] \to \Spec R$ and
$$
\sum_{k=0}^{d} dx_k \otimes X_k \in \widehat{\Omega}_{R/A} \otimes_R R[X_0,\ldots,X_d]
$$
is basic on $U(F_1,\ldots,F_r) \cap (\phi^*)^{-1}(U(\fm))$. Let $T$ denote the localization of $R[X_0,\ldots,X_d]$ at the prime ideal $\fm R[X_{0},\ldots,X_{d}]$. 
Then we have a flat local map of catenary local domains: $(R,\fm) \to (T,\fm_T)$ and $\dim R=\dim T$. By abuse of notation, we denote by $\phi^*:\Spec T \to \Spec R$ the scheme map induced by $R \to R[X_{0},\ldots,X_{d}] \to T$. Then since $\phi^*$ is faithfully flat, it is surjective.

In what follows, it suffices to consider the case when $(F_1,\ldots,F_r)T$ is a proper ideal. Let $P \subseteq T$ be any prime ideal containing $(F_1,\ldots,F_r)T$. Then it follows from Lemma \ref{prelim} together with $(\ref{zero})$ above, that $V:=V(F_1,\ldots,F_r) \subseteq \Spec T$ is a finite set and
$$
\sum_{k=0}^d dx_k \otimes X_k
$$ 
is basic on $\Spec T \backslash V$ (the complement of $\phi^*(\Spec T \backslash V)$ in $\Spec R$ is a finite set). Then it suffices to consider the proof on $U(g) \subseteq \Spec T$ for some $g \in R$ such that $\widehat{\Omega}_{R/A} \otimes_R T[g^{-1}]$ is a $T[g^{-1}]$-free module, and we reason this as follows by Noetherian induction. By applying the lemma of generic freeness to the $R/\fp$-module $\widehat{\Omega}_{R/A} \otimes_R R/\fp$ for each prime $\fp \in \Min_R(gR)$ and repeating the same discussion for $R/\fp \to T \otimes_R R/\fp$ in place of $R \to T$, we see that this process stabilizes, because $R$ is a Noetherian ring. Hence it suffices to prove the theorem on the open subset $U(g) \subseteq \Spec T$. 

Now there exists a (not necessarily homogeneous) polynomial 
$$
G \in R[X_0,\ldots,X_d] \backslash \fm R[X_0,\ldots,X_d]
$$ 
such that
$$
\sum_{k=0}^d dx_k \otimes X_k
$$
is basic on $U(g \cdot G) \subseteq \Spec R[X_0,\ldots,X_d]$. Choose a point $a=(a_0,\ldots,a_d) \in \mathbb{P}^d(A)$ together with its lift $\widetilde{a}=(\widetilde{a}_0,\ldots,\widetilde{a}_d) \in \mathbb{A}^{d+1}(A)$ such that $G(\widetilde{a}) \notin \fm$ (such a point exists, because $\# \mathbf{k}=\infty$ and $G \notin\fm R[X_0,\ldots,X_d]$). Then by Noetherian induction hypothesis, 
$$
\sum_{k=0}^d \widetilde{a}_k dx_k \in \widehat{\Omega}_{R/A}
$$
is basic on $\phi^*(\Spec T \backslash V) \subseteq \Spec R$ in view of (\cite{Fle}, Lemma 1.1). Now we get the following implication. Let $\widetilde{a} \in \mathbb{A}^{d+1}(A)$ be such that $G(\widetilde{a}) \notin \fm$. Then for all $\fp \in \phi^*(\Spec T \backslash V)=\Spec R \backslash \{\fp_1,\ldots,\fp_r,\fm\}$, it follows from Lemma \ref{lemma1} that
$$
\mathbf{x}_{\widetilde{a}}=\sum_{k=0}^{d} \widetilde{a}_k x_{k} \notin \fp^{(2)}. 
$$
So it remains to deal with the issue on a finite set $\{\fp_1,\ldots,\fp_r,\fm\}$. First, if we have $\widetilde{a}_i \in A^{\times}$ for some $i$, then $\mathbf{x}_{\widetilde{a}} \notin \fm^2$ by the minimality of $x_0,\ldots,x_d$. Therefore, we have verified the condition: $\mathbf{x}_{\widetilde{a}} \notin \fp^{(2)}$ for all $\fp \in \Spec R \backslash \{\fp_1,\ldots,\fp_r\}$, whenever $G(\widetilde{a}) \notin \fm$.

Put $\overline{U}_0:=U(\overline{G}) \subseteq \mathbb{A}^{d+1}(\mathbf{k}) \backslash\{0\}$ (the origin of $\mathbb{A}^{d+1}(\mathbf{k})$ is excluded). Then define an open subset 
\begin{equation}
\label{specificopen}
U_0 \subseteq \mathbb{P}^d(\mathbf{k})
\end{equation} 
as the image of $\overline{U}_0$ under the projection map: $\mathbb{A}^{d+1}(\mathbf{k}) \backslash\{0\} \to \mathbb{P}^d(\mathbf{k})$. To deal with the issue on $\{\fp_1,\ldots,\fp_r\}$, take the homogeneous polynomial: 
$$
F(X_0,\ldots,X_d):=\sum_{i=0}^dx_iX_i \in R[X_{0},\ldots,X_{d}].
$$
Then it suffices to force 
$$
\mathbf{x}_{\widetilde{a}}=\sum_{i=0}^d \widetilde{a}_ix_i \notin \fp_j
$$ 
for all $1 \le j \le r$. We fix a point $a=(a_0:\cdots:a_d) \in \mathbb{P}^d(A)$. For each $1 \le j \le r$, we have an equivalence of conditions:
\begin{equation}
\label{one}
\mathbf{x}_{\widetilde{a}}=\sum_{i=0}^{d} \widetilde{a}_i x_i \notin \fp_{j}~\iff F(\widetilde{a}_{0},\ldots,\widetilde{a}_{d}) \not\equiv 0 \pmod {\fp_j}.
\end{equation}
Moreover, the condition $(\ref{one})$ is stable under taking multiplication by elements of $A^{\times}$. 

Our final goal is to identify the set of points of $\mathbb{P}^d(A)$ satisfying the condition $(\ref{one})$ and describe it as the inverse image of an open subset under the map $\Sp_A:\mathbb{P}^d(A) \to \mathbb{P}^d(\mathbf{k})$. Since the maximal ideal $\fm$ of $R$ is generated by $x_0,\ldots,x_d$, it is clear that $\fm$ is generated by the set $\{F(\widetilde{a}_0,\ldots,\widetilde{a}_d)~|~(a_0 : \ldots : a_d) \in \mathbb{P}^d(A)\}$. Since the union of prime ideals $\fp_1,\ldots,\fp_r$ is strictly contained in $\fm$, there exists a point
$(a_0:\ldots :a_d) \in \mathbb{P}^d(A)$ for which $F(\widetilde{a}_0,\ldots,\widetilde{a}_d) \notin \fp_{i}$ for all $1 \le i \le r$. In fact, more is true. Applying Lemma \ref{Case1Case2} to the primes ideals $\fp_1,\ldots,\fp_r$, each of which is different from $\fm$, we find non-empty open subsets $U_1,\ldots,U_r$ of $\mathbb{P}^d(\mathbf{k})$ with each $U_i$ attached to $\fp_i$. Hence $U'$ is defined as the intersection $U_0 \cap U_1 \cap \cdots \cap U_r \subseteq \mathbb{P}^d(\mathbf{k})$, where $U_0$ is as given in $(\ref{specificopen})$. This completes the proof of the theorem.
\end{proof}

As the main result of this article, we obtain the following theorem.

\begin{theorem}[Local Bertini Theorem]\label{thm:localBertini}
\label{Bertini}
Let $(R,\fm,\mathbf{k})$ be a complete local domain of mixed characteristic $p>0$ and suppose the following conditions:

\begin{enumerate} 
\item[$\mathrm{(1)}$]
let $A \to R$ be a coefficient ring map for a complete discrete valuation ring $(A,\pi_{A})$; 

\item[$\mathrm{(2)}$]
let $x_0,x_1,\ldots,x_d$ be a fixed set of minimal generators of $\fm$;

\item[$\mathrm{(3)}$]
$R$ is normal, of $\depth R \ge 3$, and the residue field $\mathbf{k}$ is infinite.
\end{enumerate}
Then there exists a Zariski dense open subset $U \subseteq \mathbb{P}^{d}(\mathbf{k})$ satisfying the following properties. For any $a=(a_{0}:\cdots:a_{d}) \in \Sp_A^{-1}(U)$, the quotient $R/\mathbf{x}_{\widetilde{a}}R$ is a normal domain of mixed characteristic $p>0$, where we put
$$
\mathbf{x}_{\widetilde{a}}:=\sum_{i=0}^{d} \widetilde{a}_i x_i. 
$$  
\end{theorem}

\begin{proof}
The first step of the proof has been completed in Theorem \ref{theorem2}. Let $\mathrm{Reg}(R)$ denote the regular locus of $\Spec R$. Taking $U' \subseteq \mathbb{P}^d(\mathbf{k})$ as given in the conclusion of Theorem \ref{theorem2}, we have
\begin{equation}
\label{equation:ichi}
\mathrm{Reg}(R) \cap V(\mathbf{x}_{\widetilde{a}}) \subseteq \mathrm{Reg}(R/\mathbf{x}_{\widetilde{a}}R)
\end{equation}
for all $a=(a_0:\ldots:a_d) \in \Sp^{-1}_A(U')$. Let us explain basic ideas of this proof. In ${\bf{Step 1}}$, we find an open set $U'' \subseteq \mathbb{P}^d(\mathbf{k})$ to deal with finitely many bad primes, with the help of Lemma \ref{Case1Case2}. Then in ${\bf{Step 2}}$, we show that for $a \in \Sp_A^{-1}(U' \cap U'')$, the quotients $R/\mathbf{x}_{\widetilde{a}}R$ satisfy the well-known conditions $(R_1)$ and $(S_2)$ and thus are normal. Finally in ${\bf{Step 3}}$, we find an open set
$U''' \subseteq \mathbb{P}^d(\mathbf{k})$ to show that the quotient rings $R/\mathbf{x}_{\widetilde{a}}R$ with $a \in \Sp_A^{-1}(U''')$ are of mixed characteristic. 

${\bf{Step 1}}$:
The goal of this step is to find a candidate of an open subset $U'' \subseteq \mathbb{P}^d(\mathbf{k})$. Since $R$ is a complete local domain, the singular locus $\mathrm{Sing}(U(\fm))$ of $U(\fm) \subseteq \Spec R$ is a proper closed subset. Hence the set of minimal primes in $\mathrm{Sing}(U(\fm))$ is finite, and let
$$
Q_1:=\{\fp \in U(\fm)~|~\fp~\mbox{is a minimal prime in}~\mathrm{Sing}(U(\fm))\}.
$$
Note that every prime in $Q_{1}$ has height $\ge 2$, due to the $(R_1)$ condition on $R$. On the other hand,
$$
Q_2:=\{\fp \in U(\fm)~|~\mbox{$\depth R_{\fp}=2$ and $\dim R_{\fp}> 2$}\}
$$ 
is also a finite set by (\cite{Fle}, Lemma 3.2 and the $(S_2)$ condition on $R$). Now let $Q_{1} \cup Q_{2}:=\{\fp_{1},\ldots,\fp_{m}\}$ and let us put $F(X_0,\ldots,X_d)=\sum_{i=0}^{d}x_{i}X_{i} \in R[X_0,\ldots,X_d]$. Then for each $1 \le j \le m$, it follows that
$$
\mathbf{x}_{\widetilde{a}}=\sum_{i=0}^d \widetilde{a}_i x_i \notin \fp_{j}~\iff F(\widetilde{a}_0,\ldots,\widetilde{a}_d) \not \equiv 0 \pmod {\fp_j}.
$$ 
Then applying Lemma \ref{Case1Case2} to $\{\fp_1,\ldots,\fp_m\}$, we obtain $U_1,\ldots,U_m$, which are non-empty open subsets of $\mathbb{P}^d(\mathbf{k})$ with each $U_i$ attached to $\fp_i$. Put 
\begin{equation}
\label{equation:ni}
U'':=U_1 \cap \cdots \cap U_m . 
\end{equation}
 
${\bf{Step 2}}$: 
The goal of this step is to show that the quotient $R/\mathbf{x}_{\widetilde{a}}R$ is a normal domain for $a \in \Sp_A^{-1}(U'')$, or equivalently, $\mathbf{x}_{\widetilde{a}}$ is contained in no prime of $Q_{1} \cup Q_{2}$. Pick $\fp \in U(\fm) \cap V(\mathbf{x}_{\widetilde{a}})$ with $\Ht \fp \ge 2$ and assume that $\mathbf{x}_{\widetilde{a}}=\sum_{i=0}^d a_ix_i$ satisfies the condition $(\ref{equation:ichi})$ and that $\mathbf{x}_{\widetilde{a}}$ is contained in no prime of $Q_{1} \cup Q_{2}$. 

(i) If we have $\Ht \fp >2$, then since $\mathbf{x}_{\widetilde{a}}$ is contained in no prime of $Q_{2}$, it follows that $\dim (R/\mathbf{x}_{\widetilde{a}}R)_{\fp} \ge 2$ and $\depth (R/\mathbf{x}_{\widetilde{a}}R)_{\fp} \ge 2$.

(ii) If we have $\Ht \fp=2$, then since $\mathbf{x}_{\widetilde{a}}$ is contained in no prime of $Q_{1}$ and the height of every prime in $Q_{1}$ is at least 2, it follows that $R_{\fp}$ is regular. 
By $(\ref{equation:ichi})$, one finds that $(R/\mathbf{x}_{\widetilde{a}}R)_{\fp}$ is a discrete valuation ring. On the other hand, the hypothesis that $\depth (R) \ge 3$ implies that $\depth (R/\mathbf{x}_{\widetilde{a}}R) \ge 2$ by (\cite{BrHer}, Proposition 1.2.9). Hence $R/\mathbf{x}_{\widetilde{a}}R$ is a normal domain in view of Serre's normality criterion. 

${\bf{Step 3}}$:
In this final step, we make $R/\mathbf{x}_{\widetilde{a}}R$ into a local ring of mixed characteristic $p>0$.
Let $\{\fq_1,\ldots,\fq_n\}$ be a set of all height-one primes of $R$ lying above $\pi_A$. Then again, applying Lemma \ref{Case1Case2} to $\{\fq_1,\ldots,\fq_n\}$, we find a non-empty open subset 
\begin{equation}
\label{equation:san}
U''' \subseteq \mathbb{P}^d(\mathbf{k}),
\end{equation}
such that for $a=(a_0:\cdots:a_d) \in \Sp^{-1}_A(W) \subseteq \mathbf{P}^d(A)$, we have $\mathbf{x}_{\widetilde{a}} \notin \fq_1 \cup \cdots \cup \fq_n$.

Combining $(\ref{equation:ichi})$, $(\ref{equation:ni})$, together with $(\ref{equation:san})$, and taking a non-empty open subset 
$$
U:=U' \cap U'' \cap U''' \subseteq \mathbb{P}^d(\mathbf{k}),
$$
it turns out that $\Sp_A^{-1}(U) \subseteq \mathbb{P}^d(A)$ has the required property.
\end{proof}

\begin{remark}
Assume that $R$ is a Cohen-Macaulay local normal domain. Then $R/\mathbf{x}_{\widetilde{a}}R$ is Cohen-Macaulay and normal. One can continue this process until $\dim R=2$ is attained. As to Cohen-Macaulay property, the following fact is known. Assume that $S \to R$ is a torsion free module-finite extension of local domains such that $S$ is regular. Then $R$ is a flat $S$-module if and only if $R$ is Cohen-Macaulay by the Auslander-Buchsbaum formula.
\end{remark}

Next, let us consider the case when the residue field is finite. Let $(R,\fm,\mathbb{F})$ be a complete local normal domain of mixed characteristic $p>0$ with finite residue field $\mathbb{F}$. In other words, $R$ is a finite extension of $W(\mathbb{F})[[z_1,\ldots,z_n]]$, where $W(\mathbb{F})$ is the ring of Witt vectors of $\mathbb{F}$. Let $W(\mathbb{F})^{\ur}$ be the maximal unramified extension of $W(\mathbb{F})$. Then $W(\overline{\mathbb{F}})$ is the completion of $W(\mathbb{F})^{\ur}$. Put
$$
R_{W(\overline{\mathbb{F}})}:=R \widehat{\otimes}_{W(\mathbb{F})}W(\overline{\mathbb{F}})~(\mathrm{resp.}
~R_{W(\mathbb{F})^{\ur}}:=\mathrm{strict~henselization~of}~R).
$$
Then $R_{W(\mathbb{F})^{\ur}}$ is local Noetherian, but not complete. By the main result of \cite{Gre}, $R_{W(\overline{\mathbb{F}})}$ is the completion of $R_{W(\mathbb{F})^{\ur}}$ and a local normal domain. From algebraic number theory, $W(\mathbb{F})^{\ur}$ is obtained from $W(\mathbb{F})$ by adjoining all $n$-th roots of unity for $(n,p)=1$. There is a structure map $W(\overline{\mathbb{F}}) \to R_{W(\overline{\mathbb{F}})}$. We define the multiplicative map (not additive)
$$
[-]:\overline{\mathbb{F}} \to W(\overline{\mathbb{F}})
$$
as the Teichm\"uller map (see Proposition 8 in Chap. II \S 4 in \cite{Se} for details). An element in the image of $[-]$ is called a \textit{Teichm\"uller lift}. In particular, we have $q \circ [-]=\mathrm{Id}_{\overline{\mathbb{F}}}$, where $q:W(\overline{\mathbb{F}}) \to \overline{\mathbb{F}}$ is the residue field map. There is a set-theoretic mapping:
\begin{equation}
\label{Teichmullermap}
\theta_{W(\overline{\mathbb{F}})}:\mathbb{P}^d(\overline{\mathbb{F}}) \to \mathbb{P}^d(W(\overline{\mathbb{F}}))
\end{equation}
defined by $\theta_{W(\overline{\mathbb{F}})}(b):=([b_0]: \cdots: [b_d])$ for $b=(b_0:\cdots:b_d) \in 
\mathbb{P}^d(\overline{\mathbb{F}})$. Then $\theta_{W(\overline{\mathbb{F}})}$ is well-defined, since $[-]$ is multiplicative. We write 
$R_{W(\mathbb{F}')}:=R \otimes_{W(\mathbb{F})} W(\mathbb{F}')$ for a finite field extension $\mathbb{F} \to \mathbb{F}'$.

\begin{corollary}[Finite Residue Field Case]
\label{cor1}
Let the hypothesis be as in Theorem \ref{Bertini} for $(R,\fm,\mathbb{F})$, except that we now assume the residue field $\mathbb{F}$ is finite. Then there exists a non-empty open subset $U \subseteq \mathbb{P}^d(\overline{\mathbb{F}})$ such that the following holds:

Fix an element $a=(a_0:\cdots:a_d) \in \theta_{W(\overline{\mathbb{F}})}(U)$. Then there exists a finite extension $W(\mathbb{F}) \to W(\mathbb{F}')$ such that there is a choice of a lift $\widetilde{a}_i$ of $a_i$ for each $0 \le i \le d$ with $\widetilde{a}_i \in W(\mathbb{F}')$, and the quotient $R_{W(\mathbb{F}')}/\mathbf{x}_{\widetilde{a}}R_{W(\mathbb{F}')}$ is a normal domain of mixed characteristic $p>0$.
\end{corollary}

The point of the proof is to construct a multiplicative map: $\widetilde{[-]}:\overline{\mathbb{F}} \to W(\mathbb{F})^{\ur}$ which extends to the map $[-]:\overline{\mathbb{F}} \to W(\overline{\mathbb{F}})$.

\begin{proof}
We keep the notation as in Theorem \ref{Bertini}. First, note that a priori a choice of a lift $\widetilde{a}_i$ of each $a_i$ is contained in $W(\overline{\mathbb{F}})$. Since $x_0,x_1,\ldots,x_d$ are the minimal generators of the maximal ideal of $R_{W(\overline{\mathbb{F}})}$, the hypotheses of Theorem \ref{Bertini} are fulfilled for the complete local domain $R_{W(\overline{\mathbb{F}})}$. 
Let us construct a multiplicative map $\widetilde{[-]}:\overline{\mathbb{F}} \to W(\mathbb{F})^{\ur}$ which extends to the map $[-]$. Note that $W(\mathbb{F})^{\ur}=
\varinjlim_{\lambda \in \Lambda} W(\mathbb{F}_{\lambda})$ and we have the Teichm\"uller mapping $\mathbb{F}_{\lambda} \to W(\mathbb{F}_{\lambda})$. Then we have a commutative diagram:
$$
\begin{CD}
\mathbb{F}_{\lambda'} @>\mathrm{Teich}>> W(\mathbb{F}_{\lambda'}) \\
@AAA @AAA \\
\mathbb{F}_{\lambda} @>\mathrm{Teich}>>  W(\mathbb{F}_{\lambda}) \\
\end{CD}
$$
which naturally forms a direct system, so the desired map  $\widetilde{[-]}$ is given by its direct limit. On the other hand, it is easy to see that the map $[-]:\overline{\mathbb{F}} \to W(\overline{\mathbb{F}})$ factors as
$$
\begin{CD}
\overline{\mathbb{F}} @>\widetilde{[-]}>> 
W(\mathbb{F})^{\ur} @>>> W(\overline{\mathbb{F}}), 
\end{CD}
$$
and thus for any $a \in \theta_{W(\overline{\mathbb{F}})}(U)$, we have $\mathbf{x}_{\widetilde{a}}=\sum_{i=0}^d \widetilde{a}_i x_i \in R_{W(\mathbb{F})^{\ur}}$ for an appropriate choice of $\widetilde{a}_i$ for every $a_i$. Since the map
$$
R_{W(\mathbb{F})^{\ur}}/\mathbf{x}_{\widetilde{a}}R_{W(\mathbb{F})^{\ur}} \to R_{W(\overline{\mathbb{F}})}/\mathbf{x}_{\widetilde{a}}R_{W(\overline{\mathbb{F}})}
$$ 
is local flat, $R_{W(\mathbb{F})^{\ur}}/\mathbf{x}_{\widetilde{a}}R_{W(\mathbb{F})^{\ur}}$ is a local normal domain.

By what we have said above, all the coefficients of a linear form $\mathbf{x}_{\widetilde{a}}=\sum_{i=0}^d \widetilde{a}_i x_i$ are contained in some finite extension $W(\mathbb{F}) \to W(\mathbb{F}')$. In other words, for a finite \'etale extension $R \to R_{W(\mathbb{F}')}$ of normal domains, the quotient ring $R_{W(\mathbb{F}')}/\mathbf{x}_{\widetilde{a}}R_{W(\mathbb{F}')}$ is normal.
\end{proof}

\begin{remark}
Let $\phi:(R,\fm) \to (S,\fn)$ be a flat local map of local rings. Then one might think of the relationship between $R/xR$ and $S/xS$ for a non-zero divisor $x \in \fm$. In fact, in order to use the local Bertini theorem for $S$ in terms of $R$, for example, assume that $R$ and all fibers of $\phi$ are normal. Then for any $x$ such that $R/xR$ is normal, $S/xS$ is also normal. 
\end{remark}

It is important and necessary to answer the following question:

\begin{question}
\label{question}
Resume the hypothesis of Theorem \ref{Bertini} and assume that an equality holds: $\mathbf{x}_{\widetilde{a}}=u\mathbf{x}_{\widetilde{a}'}$ for a unit $u \in R^{\times}$. Then is it true that $u \in A^{\times}$?
\end{question}

This question is restated as follows: Which subset of $\mathbb{P}^d(A)$ parametrizes the set of height-one primes $\{\mathbf{x}_{\widetilde{a}}R\}$? In other words, does it give a set of distinct height-one primes of $R$? In the next section, we will answer the above question. In fact, we need to restrict our attention to the set of those points which are in the image of the Teichm\"uller mapping. This will be important in the proof of the control theorem, which will be discussed later. We end this section with an example, which applies Theorem \ref{thm:localBertini} and its proof for a given normal domain $R$.

\begin{example}
\label{algorithm}
In this example, we are dealing with the case where the residue field is finite, since we can find a non-empty set which does the job. However, the result is valid for any discrete valuation coefficient ring. Suppose that $p \ge 3$ and
$$
R:=\mathbb{Z}_p[[x_1,x_2,x_3]]/(x_1^2+x_2^2+x_3^2),
$$
which is a three-dimensional Cohen-Macaulay local normal domain. Now let us find a non-empty open subset $U \subseteq \mathbb{P}^3(\mathbb{F}_p)$ by keeping track of the proof of Theorem \ref{thm:localBertini}.

\begin{enumerate} 
\item[$\mathrm{(i)}$]
We need to have $d\mathbf{x}_{\widetilde{a}} \in \widehat{\Omega}_{R/\mathbb{Z}_p}$ basic at every $\fp \in \Spec R$. 

\item[$\mathrm{(ii)}$]
We need to determine two finite sets of primes $Q_1$ and $Q_2$ in Theorem \ref{thm:localBertini}.

\item[$\mathrm{(iii)}$]
We need to avoid $Q_1$ and $Q_2$ as above.

\end{enumerate}
Then we know
$$
\widehat{\Omega}_{R/\mathbb{Z}_p} \simeq \frac{Rdx_1\oplus Rdx_2\oplus Rdx_3}{R(x_1dx_1+x_2dx_2+x_3dx_3)},
$$
the singular locus of $R$ is defined by the ideal $(x_1,x_2,x_3)$, $Q_1=\{(x_1,x_2,x_3)\}$ and $Q_2=\varnothing$. To get a normal ring of mixed characteristic, take $\mathbf{x}_{\widetilde{a}}:=\widetilde{a}_0p+\sum_{i=1}^3 \widetilde{a}_ix_i$ such that 
$$
\overline{a}=(\overline{a}_0:\overline{a}_1:\overline{a}_2:\overline{a}_3) \in U:=U^+(z_0) \cap \big(\bigcup_{i=1}^3 U^+(z_i)\big) \subseteq \mathbb{P}^3(\mathbb{F}_p)
$$
for the homogeneous coordinate $(z_0:z_1:z_2:z_3)$. Let us check that this $U$ is what we want. We have $\mathbf{x}_{\widetilde{a}} \notin Q_1$. If we assume $\widetilde{a}_1$ is a unit for simplicity, we see that
$$
\widehat{\Omega}_{R/\mathbb{Z}_p}\Big[\frac{1}{x_3}\Big] \simeq R\Big[\frac{1}{x_3}\Big]dx_1 \oplus R\Big[\frac{1}{x_3}\Big]dx_2
$$
is a free module, in which the image of $d\mathbf{x}_{\widetilde{a}}$ spans a direct summand. On the other hand, for $\overline{R}:=R/x_3R$,
$$
\widehat{\Omega}_{R/\mathbb{Z}_p}/x_3 \cdot \widehat{\Omega}_{R/\mathbb{Z}_p} \simeq \frac{\overline{R}dx_1 \oplus \overline{R}dx_2 \oplus \overline{R}dx_3}{\overline{R}(x_1dx_1+x_2dx_2)}.
$$
To show that the image of $d\mathbf{x}_{\widetilde{a}}$ is basic on $\widehat{\Omega}_{R/\mathbb{Z}_p}/x_3 \cdot \widehat{\Omega}_{R/\mathbb{Z}_p}$, keep track of the same steps as above by inverting and killing first $x_2$, and then $x_1$.
\end{example}

\begin{remark}
If one takes $R:=\mathbb{Z}_p[[x]]$, then $p \notin \fp^{(2)}$ for every prime $\fp$ of $R$, since $p$ is a regular parameter. But $dp=0$ and so the converse of Lemma \ref{lemma1} does not hold true.
\end{remark}

\section{Distinct hyperplane sections in the local Bertini theorem}
\label{sectionsubmain}

In this section, we give an answer to Question \ref{question}. Assume that $(R,\fm,\mathbf{k})$ is a complete local (not necessarily normal) domain with perfect residue field of characteristic $p>0$ with its coefficient ring $W(\mathbf{k})$, the ring of Witt vectors. Then as defined in $(\ref{Teichmullermap})$, the mapping:
$$
\theta_{W(\mathbf{k})}:\mathbb{P}^d(\mathbf{k}) \to \mathbb{P}^d(W(\mathbf{k}))
$$ 
is induced by the Teichm\"uller mapping $\mathbf{k} \to W(\mathbf{k})$. Note that the field $\mathbf{k}$ can be finite. The following proposition asserts that the parameter set of specializations in the local Bertini theorem may be identified with a certain open subset $U \subseteq \mathbb{P}^d(\mathbf{k})$.

\begin{proposition}
\label{prop2}
Let the notation be as above and let $x_0,\ldots,x_d$ be a set of minimal generators of the maximal ideal of $R$. Assume the following conditions:
\begin{enumerate}
\item[$\mathrm{(1)}$]
If $\pi_{W(\mathbf{k})} \in \fm^2$, we put $\mathbf{x}_{\widetilde{a}}=\sum_{i=0}^d \widetilde{a}_i x_i$.

\item[$\mathrm{(2)}$]
If $\pi_{W(\mathbf{k})} \notin \fm^2$, we put $x_0=\pi_{W(\mathbf{k})}$ and  $\mathbf{x}_{\widetilde{a}}=\widetilde{a}_0\pi_{W(\mathbf{k})}+\sum_{i=1}^d \widetilde{a}_i x_i$.
\end{enumerate}
Suppose that $\mathbf{x}_{\widetilde{a}}=u\mathbf{x}_{\widetilde{a}'}$ for $a, a' \in \theta_{W(\mathbf{k})}(\mathbb{P}^d(\mathbf{k})) \subseteq \mathbb{P}^d(W(\mathbf{k}))$ and $u \in R^{\times}$. Then we have $u \in W(\mathbf{k})^{\times}$.
\end{proposition}

\begin{proof}
To clarify the notation, we simply write $a_i$ in place of $\widetilde{a}_i$.

$(1)$ Assume that $\pi_{W(\mathbf{k})} \in \fm^2$ and denote for simplicity by $a_i$ the image of $a_i \in W(\mathbf{k})$ under the surjection: 
$$
W(\mathbf{k})[[X_0,\ldots,X_d]] \twoheadrightarrow R~(X_i \mapsto x_i).
$$ 
Let $P$ be its kernel and we prove that 
\begin{equation}
\label{oneone}
P \subseteq \pi_{W(\mathbf{k})}W(\mathbf{k})[[X_0,\ldots,X_d]]+\mathcal{I}^2.
\end{equation}
Here, we put $\mathcal{I}:=(X_0,\ldots,X_d)$. Let $F \in P$ be a non-zero element. We prove $(\ref{oneone})$ by contradiction and so assume that
\begin{equation}
\label{contradiction}
F \notin \pi_{W(\mathbf{k})}W(\mathbf{k})[[X_0,\ldots,X_d]]+\mathcal{I}^2.
\end{equation}
Under this assumption, after reducing $W(\mathbf{k})[[X_0,\ldots,X_d]]$ by $\pi_{W(\mathbf{k})}$, we find that 
$$
f:=F\pmod{\pi_{W(\mathbf{k})}} \notin \mathcal{I}^2\mathbf{k}[[X_0,\ldots,X_d]].
$$ 
The number of the minimal generators of the maximal ideal of $R/\pi_{W(\mathbf{k})}R$ is equal to 
$$
\dim_{\mathbf{k}}\fm/(\pi_{W(\mathbf{k})}R+\fm^2)=\dim_{\mathbf{k}} \fm/\fm^2,
$$
because we assumed $\pi_{W(\mathbf{k})} \in \fm^2$. Hence $\overline{x}_0,\ldots,\overline{x}_d$ form the minimal generators of the maximal ideal of $R/\pi_{W(\mathbf{k})}R$.

\begin{claim}
\label{sublemma-1}
Under the assumption $(\ref{contradiction})$, we can choose $s$ with $0 \le s \le d$ such that
$$
f=\sum_{i=0}^{\infty}h_iX^i_s \in \mathbf{k}[[X_0,\ldots,X_d]],
$$ 
$h_i \in \mathbf{k}[[X_1,\ldots,X_{s-1},X_{s+1},\ldots,X_d]]$ for all $i \ge 0$, and $h_1$ is a unit.
\end{claim}

\begin{proof}[Proof of Claim \ref{sublemma-1}]
We explain how to choose $h_1$ as a unit element. Fix a presentation $f=\sum_{i=0}^{\infty}h_iX^i_s \in \mathbf{k}[[X_0,\ldots,X_d]]$ with respect to $s$ and assume that $h_1$ is not a unit. Then since $f \notin \mathcal{I}^2\mathbf{k}[[X_0,\ldots,X_d]]$, the element $h_0$ contains a non-zero linear term after presenting $h_0$ as an (infinite) sum of homogeneous polynomials. So we can write
$$
h_0=\sum_{i=1,i \ne s}^d a_i X_i+(\mbox{terms of degree at least 2})
$$
and we have $a_t \ne 0$ for some $t$. By replacing $s$ with $t$, we can achieve the requirement that $h_1$ is a unit.
\end{proof}

By mapping $f=\sum_{i=0}^{\infty}h_iX^i_s$ to the quotient $R/\pi_{W(\mathbf{k})}R$, we get 
$$
\sum_{i=1}^{\infty}\overline{h}_i\overline{x}^i_s=-\overline{h}_0~\mbox{in}~R/\pi_{W(\mathbf{k})}R. 
$$
Note that $\sum_{i=1}^{\infty}\overline{h}_i\overline{x}^{i-1}_s$ is a unit since $\overline{h}_1$ is a unit of $R/\pi_{W(\mathbf{k})}R$.
Hence, we have:
$$ 
\overline{x}_s \cdot (\mbox{unit})=-\overline{h}_0.
$$
But this gives a contradiction to the fact that $\overline{x}_0,\ldots,\overline{x}_d$ are the minimal generators of the maximal ideal of $R/\pi_{W(\mathbf{k})}R$ and that $-\overline{h}_0 \in (\overline{x}_0,\ldots,\overline{x}_{s-1},\overline{x}_{s+1},\ldots,\overline{x}_d)$. Hence $(\ref{oneone})$ is established.

In the next place, fix an arbitrary lifting $\widetilde{u} \in W(\mathbf{k})[[X_0,\ldots,X_d]]$ of $u \in R^{\times}$ under the map $W(\mathbf{k})[[X_0,\ldots,X_d]] \to R$. We write:
$$
\widetilde{u}=\sum_{s_0,\ldots,s_d}(\sum_{r=0}^{\infty}b_r^{(s_0,\ldots,s_d)} \pi_{W(\mathbf{k})}^r)X_0^{s_0}\cdots X_d^{s_d},
$$
where $(s_0,\ldots,s_d)$ denotes the multi-index and $b_r^{(s_0,\ldots,s_d)}$ are the Teichm\"uller lifts. Since $\widetilde{u}$ is a unit, we have $b_0^{(0,\ldots,0)} \ne 0$. By lifting the relation $\mathbf{x}_{\widetilde{a}}=u\mathbf{x}_{\widetilde{a}'}$ to $W(\mathbf{k})[[X_0,\ldots,X_d]]$, we have
$$
\sum_{i=0}^da_iX_i\equiv \big(\sum_{s_0,\ldots,s_d}(\sum_{r=0}^{\infty}b_r^{(s_0,\ldots,s_d)} \pi_{W(\mathbf{k})}^r)X_0^{s_0}\cdots X_d^{s_d}\big)\big(\sum_{i=0}^da'_iX_i\big)\pmod P.
$$
Rewrite the above equation as:
$$
\sum_{i=0}^d(a_i-a'_ib_0^{(0,\ldots,0)})X_i\equiv\big(\sum_{s_0,\ldots,s_d}(\sum_{r=1}^{\infty}b_r^{(s_0,\ldots,s_d)} \pi_{W(\mathbf{k})}^r)X_0^{s_0}\cdots X_d^{s_d}\big)\big(\sum_{i=0}^da'_iX_i\big)
$$
$$
+\big(\sum_{\substack{(s_0,\ldots,s_d) \\ \ne(0,\ldots,0)}}b_0^{(s_0,\ldots,s_d)}X_0^{s_0}\cdots X_d^{s_d}\big)\big(\sum_{i=0}^da'_iX_i\big) \pmod P.
$$
Then by mapping the above equation to the quotient $\mathbf{k}[[X_0,\ldots,X_d]]$, comparing the degrees on both sides, and then using the relation $(\ref{oneone})$, we find that
$$
a_i=a'_ib_0^{(0,\ldots,0)}+\pi_{W(\mathbf{k})} \cdot v_i
$$
for some $v_i \in W(\mathbf{k})$. However if $v_i \ne 0$, this implies that $a_i$ is not a Teichm\"uller lift, since both $a'_i$ and $b_0^{(0,\ldots,0)}$ are Teichm\"uller lifts. So we must have $v_i=0$ for all $0 \le i \le d$ and the following relation holds:
$$
\big(\sum_{s_0,\ldots,s_d}(\sum_{r=1}^{\infty}b_r^{(s_0,\ldots,s_d)} \pi_{W(\mathbf{k})}^r)X_0^{s_0}\cdots X_d^{s_d}\big)\big(\sum_{i=0}^da'_iX_i\big) 
$$
$$
+\big(\sum_{\substack{(s_0,\ldots,s_d) \\ \ne(0,\ldots,0)}}b_0^{(s_0,\ldots,s_d)}X_0^{s_0}\cdots X_d^{s_d}\big)\big(\sum_{i=0}^da'_iX_i\big) \in P.
$$
Since $P$ is a prime ideal, we deduce that 
$$ 
\big(\sum_{s_0,\ldots,s_d}(\sum_{r=1}^{\infty}b_r^{(s_0,\ldots,s_d)} \pi_{W(\mathbf{k})}^r)X_0^{s_0}\cdots X_d^{s_d} +\sum_{\substack{(s_0,\ldots,s_d) \\ \ne(0,\ldots,0)}}b_0^{(s_0,\ldots,s_d)}X_0^{s_0}\cdots X_d^{s_d}\big) \in P
$$ 
and thus $\widetilde{u} \equiv b_0^{(0,\ldots,0)} \pmod P$ and $u \in W(\mathbf{k})$. Now we obtain $u \in W(\mathbf{k})^{\times}$.

$(2)$ Assume that $\pi_{W(\mathbf{k})} \notin \fm^2$. Then taking $x_0=\pi_{W(\mathbf{k})}$, we may consider the surjection $W(\mathbf{k})[[X_1,\ldots,X_d]] \twoheadrightarrow R$ ($X_i \mapsto x_i$) and let $P$ be its kernel. In this case, we prove that
\begin{equation}
\label{twotwo}
P \subseteq \mathcal{M}^2,
\end{equation}
where we put $\mathcal{M}:=(\pi_{W(\mathbf{k})},X_1,\ldots,X_d)$. Let $F \in P$ be any non-zero element. Then we can write
\begin{equation}
\label{threethree}
F=\sum_{s_1,\ldots,s_d}(\sum_{r=0}^{\infty}b_r^{(s_1,\ldots,s_d)} \pi_{W(\mathbf{k})}^r)X_1^{s_1}\cdots X_d^{s_d},
\end{equation}
for Teichm\"uller lifts $b_r^{(s_1,\ldots,s_d)}$. Now assume that $F \notin \mathcal{M}^2$ for a contradiction. Then the equation $(\ref{threethree})$ can be written as
\begin{equation}
\label{fourfour}
F=b_0 \pi_{W(\mathbf{k})}+\sum_{i=1}^d b_i X_i+(\mbox{terms of degree at least 2})
\end{equation}
for Teichm\"uller lifts $b_i$ and at least one of $b_0,\ldots,b_d$ is not zero, say $b_k$. Then mapping $(\ref{fourfour})$ to $R$, $F$ goes to $0$ and we find that $x_k \in (\pi_{W(\mathbf{k})},x_2,\ldots,x_{k-1},x_k,\ldots,x_d)$, due to $b_k \ne 0$. But since $\pi_{W(\mathbf{k})},x_1,\ldots,x_d$ are the minimal generators of $\fm$, this is a contradiction. Thus, we must have $F \in \mathcal{M}^2$ and $(\ref{twotwo})$ is proved. Assume that $\mathbf{x}_{\widetilde{a}}=u\mathbf{x}_{\widetilde{a}'}$ for some $u \in R^{\times}$. Then by applying the final step of (1), together with the fact that $P \subseteq \mathcal{M}^2$, we conclude that $u \in W(\mathbf{k})^{\times}$, as desired.
\end{proof}

Now Proposition \ref{prop2} assures us that there are sufficiently many normal hyperplane sections for a local normal domain. We start with the following lemma:

\begin{lemma}
\label{Noetherian}
Let $R$ be a Noetherian domain and let $\{P_\lambda \}_{\lambda \in \Lambda}$ be an infinite set of distinct height-one primes. Then we have
$$
\bigcap_{\lambda \in \Lambda} P_\lambda=0.
$$
\end{lemma}

\begin{proof}
Assume that there is a non-zero element $a \in \bigcap_{\lambda \in \Lambda} P_\lambda$. Let $\overline{P}_\lambda$ denote the image of $P_\lambda$ under the surjection $R \twoheadrightarrow R/aR$. Then since $a \in P_\lambda$ for all $\lambda$, the set $\{\overline{P}_\lambda\}_{\lambda \in \Lambda}$ gives an infinite set of minimal prime ideals of the Noetherian ring $R/aR$. But this is a contradiction and we must have $\bigcap_{\lambda \in \Lambda} P_\lambda=0$, as desired.
\end{proof}

Applying this lemma, we get the following.

\begin{corollary}
\label{extra}
In addition to the hypothesis of Proposition \ref{prop2}, assume that $S \subseteq \mathbb{P}^d(\mathbf{k})$ is an infinite subset such that the quotient $R/\mathbf{x}_{\widetilde{a}}R$ is an integral domain for every $a \in \theta_{W(\mathbf{k})}(S)$. Then we have 
$$
\bigcap_{a \in \theta_{W(\mathbf{k})}(S)} \mathbf{x}_{\widetilde{a}}R=0.
$$
\end{corollary}

\begin{proof}
By Proposition \ref{prop2} and our assumption, we see that $\{\mathbf{x}_{\widetilde{a}}R\}_{a \in \theta_{W(\mathbf{k})}(S)}$ is an infinite set of height-one primes of a Noetherian domain $R$. Now the corollary follows from Lemma \ref{Noetherian}.
\end{proof}

It is not clear if the proof of Proposition \ref{prop2} can be modified so that it holds for any residue field. Via the proof of the proposition, we obtain the following corollary, which is useful in dealing with the unramified case.

\begin{corollary}
\label{cor2}
In addition to the hypothesis of Proposition \ref{prop2}, assume that $\pi_{W(\mathbf{k})}$ is part of minimal generators of $\fm$ and $R/\mathbf{x}_{\widetilde{a}}R$ is an integral domain of mixed characteristic for
$$
\mathbf{x}_{\widetilde{a}}=\widetilde{a}_0\pi_{W(\mathbf{k})}+\sum_{i=1}^d \widetilde{a}_i x_i
$$
and $a=(a_0:\cdots:a_d) \in \mathbb{P}^d(W(\mathbf{k}))$. Then $\pi_{W(\mathbf{k})}$ is part of minimal generators of the maximal ideal of the quotient ring $R/\mathbf{x}_{\widetilde{a}}R$.
\end{corollary}

Note that the linear form $\mathbf{x}_{\widetilde{a}}$ is not assumed to define a normal quotient $R/\mathbf{x}_{\widetilde{a}}R$.

\begin{proof}
Let $W(\mathbf{k})[[X_1,\ldots,X_d]] \twoheadrightarrow R$ be a surjective ring map with its kernel $P$ such that the images of $\pi_{W(\mathbf{k})},X_1,\ldots,X_d$ are the minimal generators of $\fm$. Then we obtain that $P \subseteq (\pi_{W(\mathbf{k})},X_1,\ldots,X_d)^2$ as in the proof of Proposition \ref{prop2}. Let 
$$
\widetilde{\mathbf{x}}_{\widetilde{a}}=\widetilde{a}_0\pi_{W(\mathbf{k})}
+\sum_{i=1}^d\widetilde{a}_iX_i \in W(\mathbf{k})[[X_1,\ldots,X_d]] 
$$ 
be a lift of $\mathbf{x}_{\widetilde{a}}$. Since $R$ and $R/\mathbf{x}_{\widetilde{a}}R$ are integral domains of mixed characteristic by assumption, $a_i \in W(\mathbf{k})$ must be a unit for some $1 \le i \le d$, which gives an isomorphism: 
$$
W(\mathbf{k})[[X_1,\ldots,X_d]]/(\widetilde{\mathbf{x}}_{\widetilde{a}}) \simeq W(\mathbf{k})[[X_1,\ldots,X_{i-1},X_{i+1},\ldots,X_d]],
$$
together with 
$$
\pi_{W(\mathbf{k})} \notin (\pi_{W(\mathbf{k})},X_1,\ldots,X_{i-1},X_{i+1},\ldots,X_d)^2.
$$
But since $\overline{P} \subseteq (\pi_{W(\mathbf{k})},X_1,\ldots,X_{i-1},X_{i+1},\ldots,X_d)^2$, where $\overline{P}$ denotes the image of $P$ in the quotient $W(\mathbf{k})[[X_1,\ldots,X_d]]/(\widetilde{\mathbf{x}}_{\widetilde{a}})$, it follows that $\pi_{W(\mathbf{k})}$ forms part of minimal generators of the maximal ideal of $R/\mathbf{x}_{\widetilde{a}}R$, as required.
\end{proof}

\begin{remark}
In this section, it is essential to assume that the residue field is perfect, which allows us to present an element in the ring $W(\mathbf{k})[[X_1,\ldots,X_d]]$ in a unique way using the Teichm\"uller mapping. Note that the results in this section as well as Theorem \ref{thm:localBertini} can be applied to a complete local normal domain $(R,\fm,\mathbf{k})$ with infinite perfect residue field and $\depth R \ge 3$. 
\end{remark}

\section{Serre's conditions $(R_n)$ and $(S_n)$}
\label{sectionmainSerre}

In this section, we prove Bertini theorems in the case when $R$ satisfies Serre's conditions on the punctured spectra of local rings. The essential part for these cases is found in the proof of Theorem \ref{thm:localBertini}. As usual, we put
$$
\mathbf{x}_{\widetilde{a}}=\sum_{i=0}^d\widetilde{a}_ix_i
$$
for $a=(a_0,\ldots,a_d) \in \mathbb{P}^d(A)$.

\begin{corollary}
\label{cor3}
Suppose that $(R,\fm,\mathbf{k})$ is a complete local reduced ring of mixed characteristic $p>0$, that conditions $(1)$ and $(2)$ of Theorem \ref{Bertini} hold, and that the residue field $\mathbf{k}$ is infinite. If $R_{\fp}$ has Serre's $(R_n)$ $($resp. $(S_n)$$)$ for all $\fp \in X$, then there exists a Zariski dense open subset $U \subseteq \mathbb{P}^{d}(\mathbf{k})$ such that for every 
$$
a=(a_{0}:\cdots :a_{d}) \in \Sp_A^{-1}(U), 
$$
the quotient $R_{\fp}/\mathbf{x}_{\widetilde{a}}R_{\fp}$ has $(R_n)$ $($resp. $(S_n)$$)$ for all $\fp \in U(\fm) \cap V(\mathbf{x}_{\widetilde{a}})$. 
\end{corollary}

\begin{proof}
We briefly sketch the proof of the corollary. Since $R$ is complete local and reduced, the non-singular locus of $R$ is non-empty. Thus, the set
$$
Q_1:=\{\fp \in U(\fm)~|~\fp~\mbox{is a minimal prime in}~\mathrm{Sing}(U(\fm))\}
$$
is finite. Let 
$$
Q_2:=\{\fp \in U(\fm)~|~\mbox{$\depth R_{\fp}=n$ and $\dim R_{\fp}>n$}\},
$$
which is also finite by (\cite{Fle}, lemma 3.2). The proof is similar to that of Theorem \ref{Bertini}. So it suffices to avoid the union of finite set of prime ideals in $Q_1 \cup Q_2$. Namely, for any $a=(a_0:\cdots:a_d) \in \Sp^{-1}_A(U)$, the localization $R_{\fp}/\mathbf{x}_{\widetilde{a}}R_{\fp}$ has $(R_n)$ (resp. $(S_n)$) for all $\fp \in U(\fm) \cap V(\mathbf{x}_{\widetilde{a}})$ and $\mathrm{Reg}(R) \cap V(\mathbf{x}_{\widetilde{a}}) \subseteq \mathrm{Reg}(R/\mathbf{x}_{\widetilde{a}}R)$.
\end{proof}

The above proof also shows that the Bertini theorem holds for mixed Serre's conditions. That is, if $R$ has $(R_s)+(S_r)$, then so does $R_{\fp}/\mathbf{x}_{\widetilde{a}}R_{\fp}$ for all $a=(a_0:\cdots:a_d) \in \Sp^{-1}_A(U)$ and all $\fp \in U(\fm) \cap V(\mathbf{x}_{\widetilde{a}})$. For instance, we obtain the Bertini theorem for reduced rings, since we know that $R$ is reduced if and only if 
$R$ has $(R_0)+(S_1)$. To be precise, we have the following version of the local Bertini theorem:

\begin{corollary}
\label{cor4}
Suppose that $(R,\fm,\mathbf{k})$ is a complete local normal domain of dimension 2 in mixed characteristic $p>0$, that conditions $(1)$ and $(2)$ of Theorem \ref{Bertini} hold, and that the residue field $\mathbf{k}$ is infinite. Then there exists a Zariski dense open subset $U \subseteq \mathbb{P}^{d}(\mathbf{k})$ such that for every 
$$
a=(a_{0}:\cdots :a_{d}) \in \Sp_A^{-1}(U),
$$
the quotient $R/\mathbf{x}_{\widetilde{a}}R$ is a reduced ring of mixed characteristic $p>0$.
\end{corollary}

\begin{proof} 
We first show that there exists a Zariski dense open subset $U' \subseteq \mathbb{P}^{d}(\mathbf{k})$ for which $R/\mathbf{x}_{\widetilde{a}}R$ is reduced for all $a \in \Sp_A^{-1}(U')$. Since $R$ is a domain, $R_{\fp}$ has $(R_0)+(S_1)$ for all $\fp \in \Spec R$. Then by Corollary \ref{cor3}, we find a desired 
$U' \subseteq \mathbb{P}^{d}(\mathbf{k})$ such that, for $a \in \Sp_A^{-1}(U')$, 
$R_{\fp}/\mathbf{x}_{\widetilde{a}}R_{\fp}$ is reduced for all $\fp \in U(\fm) \cap V(\mathbf{x}_{\widetilde{a}})$. It remains to show that the localization of $R/\mathbf{x}_{\widetilde{a}}R$ at $\fm$, which is $R/\mathbf{x}_{\widetilde{a}}R$ itself, has $(S_1)$. But since $R$ is a local normal domain of dimension 2, it is Cohen-Macaulay and thus, $R/\mathbf{x}_{\widetilde{a}}R$ has $(S_1)$. Hence $R/\mathbf{x}_{\widetilde{a}}R$ is reduced for all $a \in \Sp_A^{-1}(U')$.

To make the quotient $R/\mathbf{x}_{\widetilde{a}}R$ into a ring of mixed characteristic, we need to shrink $U'$ to an open subset $U$. To this aim, it suffices to choose $U$ such that $\pi_A$ is not a zero divisor of $R/\mathbf{x}_{\widetilde{a}}R$ for all $a \in \Sp_A^{-1}(U)$. In fact, we may take $\mathbf{x}_{\widetilde{a}}$ so that it is contained in none of prime ideals in $\mathrm{Min}_R(\pi_AR)$, since every system of parameters of $R$ is a regular sequence. Let $U''' \subseteq \mathbb{P}^{d}(\mathbf{k})$ be an open subset as in ${\bf{Step 3}}$ of Theorem \ref{Bertini}. Then the required open subset is defined as $U:=U' \cap U'''$.
\end{proof}

\begin{remark}
In place of the hypothesis of Corollary \ref{cor4}, assume that $R$ is only an integral domain. Then can one find $\mathbf{x}_{\widetilde{a}}$ such that $R/\mathbf{x}_{\widetilde{a}}R$ is also an integral domain? In the mixed characteristic case, the answer to this question is not clear yet. But there is a two-dimensional complete local normal domain over $\mathbb{C}$ without principal prime ideals at all (such an example is due to Laufer, as mentioned in \cite{Fle}. However, an explicit example is not given there). In light of this, both Corollary \ref{cor3} and Corollary \ref{cor4} seem to be the best results.
\end{remark}

\section{Characteristic ideals of torsion modules over normal domains}
\label{characteristic}

Throughout this section, we assume that $R$ is a Noetherian normal domain and $M$ is a finitely generated torsion $R$-module. Then the localization of $R$ at every height-one prime is a discrete valuation ring. We introduce an invariant of the module $M$. For an ideal $I$ of $R$, let $M[I]$ denote the maximal submodule of $M$ which is annihilated by $I$. We follow the definition of characteristic ideals by Skinner-Urban as in \cite{SkUr}. For more results and properties on characteristic ideals with its relation to the Fitting ideal, see \S~\ref{appendix}. For a finitely generated $R$-module $M$, we denote by
$$
M^{\cl}:=\Hom_R(\Hom_R(M,R),R),
$$
the \textit{reflexive closure} of $M$.

\begin{definition}\label{definition:characteristicideal}
Let the notation be as above. Then the \textit{characteristic ideal} is an ideal of $R$ defined by
$$
\Char_R(M)=\{x \in R~|~v_P(x) \ge \ell_{R_P}(M_P)~\mbox{for any height-one prime}~P\},
$$
where $v_P(-)$ is the normalized valuation of $R_P$; that is, $v_P(PR_P)=1$.
\end{definition}

Since $M$ is torsion, $\ell_{R_{\fp}}(M_{\fp})=0$ for all but finitely many height-one primes $\fp$ of $R$ and it suffices to take only height-one primes in the support of $M$ in the definition. When $M$ is not a torsion $R$-module, we put $\Char_R(M)=0$.

\begin{remark}
The formation of characteristic ideals does not commute with base change in general. For example, let $R=\mathbb{Z}_p[[x,y]]$ and let $M=R/xR$. Then $M$ is a torsion module and $\Char_R(M)=xA$. However, the $R/xR$-module $M/xM$ is not torsion. Therefore, 
$$
0=\Char_{R/xR}(M/xM) \subsetneq xR=\Char_R(M)R/xR
$$ 
in this case. In general, even when $M/xM$ is a torsion $R/xR$-module, it may happen that $\Char_{R/xR}(M/xM) \ne \Char_R(M)R/xR$, which is caused by the presence of pseudo-null submodules (see the definition below).
\end{remark}

\begin{definition}
A finitely generated module $M$ over a Noetherian normal domain $R$ is \textit{pseudo-null}, if $M_{\fp}=0$ for all height-one primes $\fp \in \Spec R$. A homomorphism of $R$-modules $f:M \to N$ is a \textit{pseudo-isomorphism}, if both $\ker(f)$ and $\coker(f)$ are pseudo-null modules.
\end{definition}

The proof of the next lemma is found in (\cite{NSW}, Proposition 5.1.7).

\begin{lemma}[Structure Theorem]
\label{lemma5}
Let $M$ be a finitely generated torsion module over a Noetherian normal domain $R$. Then there exist a finite set of height-one primes $\{P_i\}_{i \in I}$ $($which is not necessarily a redundant set of height-one primes$)$ and a set of natural numbers $\{e_i\}_{i \in I}$ such that there is a homomorphism: 
$$
f:M \to \bigoplus_{i \in I}R/P_i^{e_i}
$$ 
that is a pseudo-isomorphism. Moreover, both $\{P_i\}_{i \in I}$ and $\{e_i\}_{i \in I}$ are uniquely determined. 
\end{lemma}

Henceforth, we use the notation $M \approx N$ to indicate that there is a pseudo-isomorphism between $M$ and $N$. As mentioned before, the formation of characteristic ideals does not commute with base change in general, which can produce extra error terms.

\begin{proposition}
\label{prop3}
Let $M$ be a finitely generated torsion module over a Noetherian normal domain $R$. Let $x$ be an element of $R$ which satisfies the following conditions: 
\begin{enumerate}
\item[$\mathrm{(1)}$]
$R/xR$ is a normal domain $($which implies that $xR$ is a prime ideal$)$; 

\item[$\mathrm{(2)}$]
$x$ is contained in no prime ideal of height-one in the support of $M$ $($which implies that $M/xM$ is a torsion $R/xR$-module$)$. 
\end{enumerate}
Then we have:
$$
\Char_{R/xR}(M/xM)=\Big(\Char_{R/xR}(M[x])\cdot\prod_{\substack{\Ht\fp=1,\\ \fp \in \Spec R}} (\fp(R/xR))^{\ell_{R_{\fp}}(M_{\fp})}\Big)^{\cl}.
$$
\end{proposition}

\begin{proof}
Recall from Proposition \ref{prop5} in Appendix, that if $0 \to L \to M \to N \to 0$ is a short exact sequence of finitely generated $R$-modules, then $\Char_R(M)=(\Char_R(L) \cdot \Char_R(N))^{\cl}$. The condition stated in the proposition implies that all relevant modules are torsion. First off, we claim that
\begin{equation}
\label{equality}
\Char_{R/xR}\Big(\bigoplus_{\Ht\fp=1} R/(\fp+xR)^{\ell_{R_{\fp}}(M_{\fp})}\Big)=\Big(\prod_{\substack{\Ht\fp=1,\\ \fp \in \Spec R}} (\fp(R/xR))^{\ell_{R_{\fp}}(M_{\fp})}\Big)^{\cl}.
\end{equation}
For the proof of $(\ref{equality})$, let $\{P_1,\ldots,P_m\}$ be a set of all height-one primes of $R/xR$ which contain $\fp(R/xR)$. Then the localization $(R/xR)_{P_i}$ is a discrete valuation ring and we have
$$
\fp(R/xR)_{P_i}=P_i^{\ell_{(R/xR)_{P_i}}((R/(\fp+xR))_{P_i})}(R/xR)_{P_i}
$$
for all $i$. Then $(\ref{equality})$ follows from this and the definition of characteristic ideals. Thus, it suffices to prove the following equality:
{\small \begin{equation}
\label{equality2}
\Big(\Char_{R/xR}(M/xM) \cdot \big(\Char_{R/xR}(M[x])\big)^{-1}\Big)^{\cl}=\Char_{R/xR}\Big(\bigoplus_{\substack{\Ht\fp=1,\\ \fp \in \Spec R}} R/(\fp+xR)^{\ell_{R_{\fp}}(M_{\fp})}\Big),
\end{equation}}
which is regarded as an element in the group of reflexive fractional ideals of $R$, and $\big(\Char_{R/xR}(M[x])\big)^{-1}$ is the unique fractional ideal which is an inverse of $\Char_{R/xR}(M[x])$.

By considering all torsion $R$-modules satisfying the condition $(2)$ with respect to $x \in R$, we show that both sides of $(\ref{equality2})$, regarded as operations on such $R$-modules, is multiplicative on short exact sequences. Let
$$
0 \to L \to M \to N \to 0
$$
be a short exact sequence of $R$-modules satisfying the condition $(2)$. Then since the function $\ell(-)$ is additive, it follows that the right hand side of $(\ref{equality2})$ is multiplicative. On the other hand, there follows the exact sequence:
$$
0 \to L[x] \to M[x] \to N[x] \to L/xL \to M/xM \to N/xN \to 0
$$
by the snake lemma. If $\fp$ is a height-two prime ideal of $R$ containing $x \in R$, we may localize the above exact sequence at $\fp$, so it follows that the left hand side of $(\ref{equality2})$ is multiplicative as well.

Hence we are reduced to the case that $M=R/\fq$ for a prime ideal $\fq$ by the prime filtration argument. Assume first that $\Ht\fq=1$. Then we clearly have $\Supp(M)=\{\fq\}$, $M[x]=0$, and $\ell_{R_\fq}(M_{\fq})=1$, because $M_{\fq}$ is a simple $R$-module. Now $(\ref{equality2})$ is obviously true. Next assume that $\Ht \fq >1$. Then it is easy to see that both sides of $(\ref{equality2})$ are equal to a unit ideal, which completes the proof.
\end{proof}

\section{Applications to characteristic ideals}
\label{application}

Our final goal is to prove Theorem \ref{charid}. The aim of the main theorem is to establish some techniques, which enable us to study the Iwasawa's main conjecture attached to a $p$-adic family of modular forms (see a forthcoming paper \cite{OcSh}). For this reason, it is necessary to deal with local rings with finite residue field. Let $(R,\fm,\mathbb{F})$ be a complete local normal domain of mixed characteristic with finite residue field. In other words, $R$ is the integral closure of $\mathbb{Z}_p[[z_1,\ldots,z_n]]$ in a finite field extension of the field of fractions of $\mathbb{Z}_p[[z_1,\ldots,z_n]]$. We recall the set-up of Corollary \ref{cor1} and prove some preliminary results.

Let $R_{W(\overline{\mathbb{F}})}:=R \widehat{\otimes}_{W(\mathbb{F})}W(\overline{\mathbb{F}})$ with its coefficient ring  $W(\overline{\mathbb{F}})$. Then if $\depth R \ge 3$, the complete local ring $R_{W(\overline{\mathbb{F}})}$ fits into the hypothesis of Theorem \ref{Bertini}. We have the set-theoretic mapping: $\theta_{W(\overline{\mathbb{F}})}:\mathbb{P}^d(\overline{\mathbb{F}}) \to \mathbb{P}^d(W(\overline{\mathbb{F}}))$ as constructed in $(\ref{Teichmullermap})$.

To establish the fundamental theorem for characteristic ideals, we need to relate torsion $R$-modules to torsion $R_{W(\overline{\mathbb{F}})}$-modules and then descend to $R$ by faithful flatness. The advantage of working with $R_{W(\overline{\mathbb{F}})}$ is that the residue field is the algebraic closure of a finite field. We introduce some notation. Denote by $\Fitt_A(M)$ the Fitting ideal of an $A$-module $M$. We make free use of results and notation from Appendix.

\begin{definition}
\label{definition:L}
Under the notation as above, fix a set of minimal generators $x_0,\ldots,x_d$ of the unique maximal ideal of $R$ and let $U \subseteq \mathbb{P}^d(\overline{\mathbb{F}})$ be as given in Corollary \ref{cor1}. Then we set
$$
\mathcal{L}_{W(\overline{\mathbb{F}})}:=\{\mathbf{x}_{\widetilde{a}}R_{W(\overline{\mathbb{F}})}~|~a=(a_0:\cdots:a_d) \in \theta_{W(\overline{\mathbb{F}})}(U)\}
$$
for the mapping $\theta_{W(\overline{\mathbb{F}})}:\mathbb{P}^d(\overline{\mathbb{F}}) \to \mathbb{P}^d(W(\overline{\mathbb{F}}))$. For a finitely generated torsion 
$R_{W(\overline{\mathbb{F}})}$-module $M$, we define a subset: 
$$
\mathcal{L}_{W(\overline{\mathbb{F}})}(M) \subseteq \mathcal{L}_{W(\overline{\mathbb{F}})}
$$ 
which consists of all height-one primes $\mathbf{x}_{\widetilde{a}}R_{W(\overline{\mathbb{F}})} \in \mathcal{L}_{W(\overline{\mathbb{F}})}$ such that the following conditions are satisfied:

\begin{enumerate}
\item[\textbf{(A)}]
$M/\mathbf{x}_{\widetilde{a}}M$ is a torsion $R_{W(\overline{\mathbb{F}})}/\mathbf{x}_{\widetilde{a}}R_{W(\overline{\mathbb{F}})}$-module.

\item[\textbf{(B)}]
The following equalities of ideals hold in $R_{W(\overline{\mathbb{F}})}/\mathbf{x}_{\widetilde{a}}R_{W(\overline{\mathbb{F}})}$:
$$
\Char_{R_{W(\overline{\mathbb{F}})}/\mathbf{x}_{\widetilde{a}}R_{W(\overline{\mathbb{F}})}}(M/\mathbf{x}_{\widetilde{a}}M)=\Big(\Char_{R_{W(\overline{\mathbb{F}})}}(M)(R_{W(\overline{\mathbb{F}})}/\mathbf{x}_{\widetilde{a}}R_{W(\overline{\mathbb{F}})})\Big)^{\cl}
$$
$$
=\Big(\Fitt_{R_{W(\overline{\mathbb{F}})}}(\bigoplus_{i=1}^m R_{W(\overline{\mathbb{F}})}/P_i^{e_i})(R_{W(\overline{\mathbb{F}})}/
\mathbf{x}_{\widetilde{a}}R_{W(\overline{\mathbb{F}})})\Big)^{\cl},
$$
where $M \approx \bigoplus_{i=1}^m R_{W(\overline{\mathbb{F}})}/P_i^{e_i}$ is a fundamental pseudo-isomorphism.
\end{enumerate}
\end{definition}

When the base ring $R$ is isomorphic to a power-series ring over a discrete valuation ring, compare the above definition with the one given in (\cite{Oc1}, Definition 3.2.). The way of interpreting the condition \textbf{(B)} in Definition \ref{definition:L} is that one wants to consider the characteristic ideals through Fitting ideal as an intermediate invariant. Note that all of three ideals appearing in \textbf{(B)} may differ in general. 

We are going to prove a number of lemmas to describe an explicit structure of the set 
$\mathcal{L}_{W(\overline{\mathbb{F}})}(M)$ (see Lemma \ref{lemma7} and Lemma \ref{lemma8}).
To this aim, we prove some preliminary results from \cite{Oc1} over general normal domains by making necessary modifications. The next lemma is necessary for a technical reason and it will be shown in Lemma \ref{Fitting5} how to use it in a more concrete situation. The reader may skip the proof at the first reading.

\begin{lemma}
\label{lemma6}
Under the notation and the hypothesis as in Corollary \ref{cor1}, assume that $M$ is a finitely generated torsion $R_{W(\overline{\mathbb{F}})}$-module. For a fundamental pseudo-isomorphism $M \approx \bigoplus_{i=1}^m R_{W(\overline{\mathbb{F}})}/P_i^{e_i}$, set
$$
I=\Fitt_{R_{W(\overline{\mathbb{F}})}}(\bigoplus_{i=1}^m R_{W(\overline{\mathbb{F}})}/P_i^{e_i})
$$
and consider the natural injection $I \to I^{\cl}$ with its cokernel $N$. Then there exists a finite set $\{Q_i\}_{1 \le i \le \ell}$ consisting of height-two primes of $R_{W(\overline{\mathbb{F}})}$ with the following condition:

Fix an arbitrary element
$$
\mathbf{x}_{\widetilde{a}}R_{W(\overline{\mathbb{F}})} \in \bigcap_{1 \le i \le \ell} \mathcal{L}_{W(\overline{\mathbb{F}})}(R_{W(\overline{\mathbb{F}})}/Q_i)
$$
and let $P \in \Spec R_{W(\overline{\mathbb{F}})}$ be any prime ideal such that $\mathbf{x}_{\widetilde{a}} \in P$ and $\Ht P \le 2$. Then $N_P=0$.
\end{lemma}

\begin{proof}
By definition, the module $N$ is supported on a closed subset of codimension two in $\Spec R_{W(\overline{\mathbb{F}})}$. So there are only finitely many height-two primes contained in $\Supp N$. Then we show that it is sufficient to choose $Q_1,\ldots,Q_{\ell}$ as those height-two primes contained in $\Supp N$. Let $N_i:=R_{W(\overline{\mathbb{F}})}/Q_i$. Then $N_i$ is a pseudo-null $R_{W(\overline{\mathbb{F}})}$-module and we have $\Char_{R_{W(\overline{\mathbb{F}})}}(N_i)=R_{W(\overline{\mathbb{F}})}$. Then by the condition \textbf{(B)}, we have 
$$
\Char_{R_{W(\overline{\mathbb{F}})}/\mathbf{x}_{\widetilde{a}}R_{W(\overline{\mathbb{F}})}}(N_i/\mathbf{x}_{\widetilde{a}}N_i)=R_{W(\overline{\mathbb{F}})}/\mathbf{x}_{\widetilde{a}}R_{W(\overline{\mathbb{F}})},
$$
and this implies that $N_i/\mathbf{x}_{\widetilde{a}}N_i$ is a pseudo-null $R_{W(\overline{\mathbb{F}})}/\mathbf{x}_{\widetilde{a}}
R_{W(\overline{\mathbb{F}})}$-module. In other words, $\mathbf{x}_{\widetilde{a}} \notin Q_i$ for all $1 \le i \le \ell$. So if we choose a prime ideal $P$ such that $\mathbf{x}_{\widetilde{a}} \in P$ and $\Ht P \le 2$, then we must have $P \notin \Supp N$, which proves the lemma.
\end{proof}

We will discuss when equalities occur between various ideals in the condition \textbf{(B)} in Definition \ref{definition:L}.

\begin{discussion}
\label{Fitting}
Suppose $M$ is a finitely generated torsion module over a Noetherian normal domain $A$. Then we have
$$
\Fitt_A(M) \subseteq \Char_A(M)~\mbox{and}~\Fitt_A(M)^{\cl}=\Char_A(M)
$$ 
(see Proposition \ref{prop5} in Appendix). Take a fundamental pseudo-isomorphism:
$$
M \approx \bigoplus_{i=1}^m R_{W(\overline{\mathbb{F}})}/P_i^{e_i}
$$
for a (not necessarily redundant) finite set of height-one primes $\{P_i\}$ of $R_{W(\overline{\mathbb{F}})}$. Choose $\mathbf{x}_{\widetilde{a}} \in R_{W(\overline{\mathbb{F}})}$ such that $R_{W(\overline{\mathbb{F}})}/\mathbf{x}_{\widetilde{a}}R_{W(\overline{\mathbb{F}})}$ is normal and $M/\mathbf{x}_{\widetilde{a}}M$ is a torsion $R_{W(\overline{\mathbb{F}})}/\mathbf{x}_{\widetilde{a}}R_{W(\overline{\mathbb{F}})}$-module. In particular, the multiplication map:
$$
\begin{CD}
\bigoplus_{i=1}^m R_{W(\overline{\mathbb{F}})}/P_i^{e_i} @>\mathbf{x}_{\widetilde{a}}>> \bigoplus_{i=1}^m R_{W(\overline{\mathbb{F}})}/P_i^{e_i}\\
\end{CD}
$$
is injective. Then since $\Fitt_B(M \otimes_AB)=\Fitt_A(M)B$ for any Noetherian $A$-algebra $B$, letting $A=R_{W(\overline{\mathbb{F}})}$ and $B=R_{W(\overline{\mathbb{F}})}/\mathbf{x}_{\widetilde{a}}R_{W(\overline{\mathbb{F}})}$, we have
\begin{equation}\label{Fitt1}
\text{{\footnotesize $ 
\Big(\Fitt_{R_{W(\overline{\mathbb{F}})}}(\bigoplus_{i=1}^m R_{W(\overline{\mathbb{F}})}/P_i^{e_i})(R_{W(\overline{\mathbb{F}})}/\mathbf{x}_{\widetilde{a}}R_{W(\overline{\mathbb{F}})})\Big)^{\cl}=\Big(\Fitt_{R_{W(\overline{\mathbb{F}})}/\mathbf{x}_{\widetilde{a}}R_{W(\overline{\mathbb{F}})}}(\bigoplus_{i=1}^m R_{W(\overline{\mathbb{F}})}/(\mathbf{x}_{\widetilde{a}}R_{W(\overline{\mathbb{F}})}+P_i^{e_i}))\Big)^{\cl}.
$}}
\end{equation}
Let $I:=\Fitt_{R_{W(\overline{\mathbb{F}})}}(\bigoplus_{i=1}^m R_{W(\overline{\mathbb{F}})}/P_i^{e_i})$. Then we have a commutative square:
$$
\begin{CD}
I(R_{W(\overline{\mathbb{F}})}/\mathbf{x}_{\widetilde{a}}R_{W(\overline{\mathbb{F}})}) @>\phi>> I^{\cl}(R_{W(\overline{\mathbb{F}})}/\mathbf{x}_{\widetilde{a}}R_{W(\overline{\mathbb{F}})}) \\
@| @| \\
\Fitt_{R_{W(\overline{\mathbb{F}})}}(\bigoplus_{i=1}^m R_{W(\overline{\mathbb{F}})}/P_i^{e_i})(R_{W(\overline{\mathbb{F}})}/\mathbf{x}_{\widetilde{a}}R_{W(\overline{\mathbb{F}})}) @>>> \Char_{R_{W(\overline{\mathbb{F}})}}(M)(R_{W(\overline{\mathbb{F}})}/\mathbf{x}_{\widetilde{a}}R_{W(\overline{\mathbb{F}})}) \\
\end{CD}
$$
Applying Lemma \ref{lemma6} to the inclusion $I \to I^{\cl}$, we have a set of height-two primes $\{Q_i\}_{1 \le i \le \ell}$. Therefore, the localization of $\phi$ at any height-one prime of $R_{W(\overline{\mathbb{F}})}/\mathbf{x}_{\widetilde{a}}R_{W(\overline{\mathbb{F}})}$ is an isomorphism if and only if the following condition holds:
\begin{equation}
\label{Fitt2}
\mathbf{x}_{\widetilde{a}}R_{W(\overline{\mathbb{F}})} \in \bigcap_{1 \le i \le \ell} \mathcal{L}_{W(\overline{\mathbb{F}})}(R_{W(\overline{\mathbb{F}})}/Q_i).
\end{equation}
We conclude that if $(\ref{Fitt2})$ is satisfied, then we have
\begin{equation}
\label{Fitt3}
\text{{\small $
\Big(\Char_{R_{W(\overline{\mathbb{F}})}}(M)(R_{W(\overline{\mathbb{F}})}/\mathbf{x}_{\widetilde{a}}R_{W(\overline{\mathbb{F}})})\Big)^{\cl}
=\Big(\Fitt_{R_{W(\overline{\mathbb{F}})}}(\bigoplus_{i=1}^m R_{W(\overline{\mathbb{F}})}/P_i^{e_i})(R_{W(\overline{\mathbb{F}})}
/\mathbf{x}_{\widetilde{a}}R_{W(\overline{\mathbb{F}})})\Big)^{\cl}.
$}}
\end{equation}
It also yields that
 
\begin{equation}
\label{Fitt4}
\text{{\small $
\Big(\Char_{R_{W(\overline{\mathbb{F}})}}(M)(R_{W(\overline{\mathbb{F}})}/\mathbf{x}_{\widetilde{a}}R_{W(\overline{\mathbb{F}})})\Big)^{\cl}=\Char_{R_{W(\overline{\mathbb{F}})}/\mathbf{x}_{\widetilde{a}}R_{W(\overline{\mathbb{F}})}}\Big(\bigoplus_{i=1}^m R_{W(\overline{\mathbb{F}})}/(\mathbf{x}_{\widetilde{a}}R_{W(\overline{\mathbb{F}})}+P_i^{e_i})\Big).
$}}
\end{equation}
\end{discussion}

We state the conclusion of the above discussion as a lemma.

\begin{lemma}
\label{Fitting5}
Let the set-up be as in Discussion \ref{Fitting}. Then the equation $(\ref{Fitt1})$ holds true without any condition and if $(\ref{Fitt2})$ is satisfied, then the equations $(\ref{Fitt3})$ and $(\ref{Fitt4})$ hold true.
\end{lemma}

We prove a lemma which describes the set $\mathcal{L}_{W(\overline{\mathbb{F}})}(M)$.

\begin{lemma}
\label{lemma7}
Under the notation and the hypothesis as in Corollary \ref{cor1}, assume that $M$ is a finitely generated torsion $R_{W(\overline{\mathbb{F}})}$-module. Then the following assertions hold:

\begin{enumerate}
\item[$\mathrm{(1)}$]
The set $\mathcal{L}_{W(\overline{\mathbb{F}})}(M)$ is identified with the intersection:
\small 
$$
\{\mathbf{x}_{\widetilde{a}}R_{W(\overline{\mathbb{F}})}~|~M/\mathbf{x}_{\widetilde{a}}M~\mathrm{is~a~torsion}~R_{W(\overline{\mathbb{F}})}/\mathbf{x}_{\widetilde{a}}
R_{W(\overline{\mathbb{F}})}\mbox{-}\mathrm{module}\} \cap \mathcal{L}_{W(\overline{\mathbb{F}})}(M_{\mathrm{null}}) \cap \bigcap_{1 \le i \le \ell}\mathcal{L}_{W(\overline{\mathbb{F}})}(R_{W(\overline{\mathbb{F}})}/Q_i),
$$
\normalsize 
where $M_{\mathrm{null}}$ is the maximal pseudo-null submodule of $M$ and $\{Q_i\}_{1 \le i \le \ell}$ is a set of height-two primes attached to $M$ as stated in Lemma \ref{lemma6}.

\item[$\mathrm{(2)}$]
Assume that $N$ is a finitely generated pseudo-null $R_{W(\overline{\mathbb{F}})}$-module and $\{Q'_i\}_{1 \le i \le k}$ is a set of all associated prime ideals of height two for the module $N$. Then we have:
$$
\mathcal{L}_{W(\overline{\mathbb{F}})}(N)=\bigcap_{1 \le i \le k} \mathcal{L}_{W(\overline{\mathbb{F}})}(R_{W(\overline{\mathbb{F}})}/Q'_i).
$$ 
\end{enumerate}
\end{lemma}

\begin{proof}
$(1)$ This is taken from (\cite{Oc1}, Lemma 3.4.), but we give its proof, as it requires some modifications. Let $M_{\mathrm{null}}$ be the maximal pseudo-null submodule of $M$. Then by Lemma \ref{lemma5}, there is a fundamental pseudo-isomorphism $M \to \bigoplus_{i=1}^m R_{W(\overline{\mathbb{F}})}/P_i^{e_i}$, together with the following commutative diagram with exact rows:
$$
\begin{CD}
0 @>>> M/M_{\mathrm{null}} @>>> \bigoplus_{i=1}^m R_{W(\overline{\mathbb{F}})}/P_i^{e_i} @>>> N @>>> 0\\
@. @V\mathbf{x}_{\widetilde{a}}VV @V\mathbf{x}_{\widetilde{a}}VV @V\mathbf{x}_{\widetilde{a}}VV \\
0 @>>> M/M_{\mathrm{null}} @>>> \bigoplus_{i=1}^m R_{W(\overline{\mathbb{F}})}/P_i^{e_i} @>>> N @>>> 0\\
\end{CD}
$$
for a (not necessarily redundant) set of heigh-one primes $\{P_i\}_{1 \le i \le m}$ of $R_{W(\overline{\mathbb{F}})}$ and $N$ is a pseudo-null module. 

First, we note that $\mathcal{L}_{W(\overline{\mathbb{F}})}(M)$ is contained in the intersection as stated in the lemma. So we need to establish the other inclusion. Now assume that
$\mathbf{x}_{\widetilde{a}}R_{W(\overline{\mathbb{F}})}$ satisfies the following conditions:
\begin{enumerate}
\item[$\bullet$]
$M/\mathbf{x}_{\widetilde{a}}M$ is a torsion $R_{W(\overline{\mathbb{F}})}/\mathbf{x}_{\widetilde{a}}R_{W(\overline{\mathbb{F}})}$-module.

\item[$\bullet$]
$\mathbf{x}_{\widetilde{a}}R_{W(\overline{\mathbb{F}})} \in \bigcap_{1 \le i \le \ell} \mathcal{L}_{W(\overline{\mathbb{F}})}(R_{W(\overline{\mathbb{F}})}/Q_i)$.
\end{enumerate}
Then the set of all elements  $\mathbf{x}_{\widetilde{a}}R_{W(\overline{\mathbb{F}})}$ satisfying the above conditions contains $\mathcal{L}_{W(\overline{\mathbb{F}})}(M)$ as a subset. So assuming that $\mathbf{x}_{\widetilde{a}}R_{W(\overline{\mathbb{F}})}$ satisfies the above conditions, it suffices to prove the following implication:
\begin{equation}
\label{seq0}
\mathbf{x}_{\widetilde{a}}R_{W(\overline{\mathbb{F}})} \in \mathcal{L}_{W(\overline{\mathbb{F}})}(M_{\mathrm{null}})
\Rightarrow
\mathbf{x}_{\widetilde{a}}R_{W(\overline{\mathbb{F}})} \in \mathcal{L}_{W(\overline{\mathbb{F}})}(M).
\end{equation}
We establish $(\ref{seq0})$ below. By our choice, the multiplication map
$$
\bigoplus_{i=1}^m R_{W(\overline{\mathbb{F}})}/P_i^{e_i} \xrightarrow{\mathbf{x}_{\widetilde{a}}} \bigoplus_{i=1}^m R_{W(\overline{\mathbb{F}})}/P_i^{e_i}\\
$$
is injective. So the snake lemma yields the following exact sequence:
\begin{equation}
\label{seq1}
0 \to N[\mathbf{x}_{\widetilde{a}}] \to M/(\mathbf{x}_{\widetilde{a}}M+M_{\mathrm{null}}) \to \bigoplus_{i=1}^m R_{W(\overline{\mathbb{F}})}/(\mathbf{x}_{\widetilde{a}}R_{W(\overline{\mathbb{F}})}+P_i^{e_i}) \to N/\mathbf{x}_{\widetilde{a}}N \to 0.
\end{equation}
There is a short exact sequence:
\begin{equation}
\label{seq2}
0 \to N[\mathbf{x}_{\widetilde{a}}] \to N \xrightarrow{\mathbf{x}_{\widetilde{a}}} N \to N/\mathbf{x}_{\widetilde{a}}N \to 0 \\
\end{equation}
of pseudo-null $R_{W(\overline{\mathbb{F}})}$-modules, where both $N[\mathbf{x}_{\widetilde{a}}]$ and $N/\mathbf{x}_{\widetilde{a}}N$ are naturally $R_{W(\overline{\mathbb{F}})}/\mathbf{x}_{\widetilde{a}}R_{W(\overline{\mathbb{F}})}$-modules. Localizing both $(\ref{seq1})$ and $(\ref{seq2})$ at all height-two primes $P$ of $R_{W(\overline{\mathbb{F}})}$ containing $\mathbf{x}_{\widetilde{a}}$, a length computation for the sequence localized at $P$ reveals that
\begin{equation}
\label{seq3}
\text{{\footnotesize $
\Char_{R_{W(\overline{\mathbb{F}})}/\mathbf{x}_{\widetilde{a}}R_{W(\overline{\mathbb{F}})}}\Big(M/(\mathbf{x}_{\widetilde{a}}M+M_{\mathrm{null}})\Big)=\Char_{R_{W(\overline{\mathbb{F}})}/\mathbf{x}_{\widetilde{a}}R_{W(\overline{\mathbb{F}})}} \Big(\bigoplus_{1 \le i \le m} R_{W(\overline{\mathbb{F}})}/(\mathbf{x}_{\widetilde{a}}R_{W(\overline{\mathbb{F}})}+P_i^{e_i})\Big).
$}}
\end{equation}
On the other hand, Lemma \ref{Fitting5} shows that
\begin{equation}
\label{seq4}
\text{{\small $
\Big(\Char_{R_{W(\overline{\mathbb{F}})}}(M)(R_{W(\overline{\mathbb{F}})}/\mathbf{x}_{\widetilde{a}}R_{W(\overline{\mathbb{F}})})\Big)^{\cl}
=\Char_{R_{W(\overline{\mathbb{F}})}/\mathbf{x}_{\widetilde{a}}R_{W(\overline{\mathbb{F}})}} \Big(\bigoplus_{1 \le i \le m} R_{W(\overline{\mathbb{F}})}/(\mathbf{x}_{\widetilde{a}}R_{W(\overline{\mathbb{F}})}+P_i^{e_i})\Big).
$}}
\end{equation}
Then combining both $(\ref{seq3})$ and $(\ref{seq4})$, we get
\begin{equation}
\label{seq5}
\Big(\Char_{R_{W(\overline{\mathbb{F}})}}(M)(R_{W(\overline{\mathbb{F}})}/
\mathbf{x}_{\widetilde{a}}R_{W(\overline{\mathbb{F}})})\Big)^{\cl}=
\Char_{R_{W(\overline{\mathbb{F}})}/\mathbf{x}_{\widetilde{a}}R_{W(\overline{\mathbb{F}})}}
\Big(M/(\mathbf{x}_{\widetilde{a}}M+M_{\mathrm{null}})\Big)
\end{equation}

Finally, since the multiplication on $M/M_{\mathrm{null}}$ by $\mathbf{x}_{\widetilde{a}}$ is injective, we have an exact sequence:
\begin{equation}
\label{seq6}
0\to M_{\mathrm{null}}/\mathbf{x}_{\widetilde{a}}M_{\mathrm{null}} \to M/\mathbf{x}_{\widetilde{a}}M \to M/(\mathbf{x}_{\widetilde{a}}M+M_{\mathrm{null}}) \to 0.
\end{equation}
Taking characteristic ideals to $(\ref{seq6})$, in view of $(\ref{seq5})$ and the condition \textbf{(B)} in Definition \ref{definition:L}, it follows that the desired implication $(\ref{seq0})$ holds if $M_{\mathrm{null}}/\mathbf{x}_{\widetilde{a}}M_{\mathrm{null}}$ is
a pseudo-null $R_{W(\overline{\mathbb{F}})}
\mathbf{x}_{\widetilde{a}}R_{W(\overline{\mathbb{F}})}$-module, which is true if $\mathbf{x}_{\widetilde{a}}R_{W(\overline{\mathbb{F}})} \in \mathcal{L}_{W(\overline{\mathbb{F}})}(M_{\mathrm{null}})$.
 This completes the proof of $(1)$.

$(2)$ This part is done in (\cite{Oc1}, Lemma 3.5 together with Lemma 3.1) in case $R_{W(\overline{\mathbb{F}})}$ is regular, so we leave the proof with necessary modifications to the reader.
\end{proof}

The next purpose is to show that $\mathcal{L}_{W(\overline{\mathbb{F}})}(M)$ is an infinite set for a finitely generated torsion $R_{W(\overline{\mathbb{F}})}$-module $M$. We need to have sufficiently many specializations of $R_{W(\overline{\mathbb{F}})}$ that are normal in order to prove the control theorem by combining Proposition \ref{prop2}, Lemma \ref{lemma6}, and Lemma \ref{lemma7}.

\begin{lemma} 
\label{lemma8}
Under the notation and the hypothesis as in Corollary \ref{cor1}, assume that $M$ is a finitely generated torsion $R_{W(\overline{\mathbb{F}})}$-module. Then we have the following assertions:

\begin{enumerate}
\item[$\mathrm{(1)}$]
The subset $\mathcal{L}_{W(\overline{\mathbb{F}})}(M) \subseteq \theta_{W(\overline{\mathbb{F}})}(U)$ may be identified with a non-empty open subset of $\mathbb{P}^d(\overline{\mathbb{F}})$ under the mapping $\theta_{W(\overline{\mathbb{F}})}:\mathbb{P}^d(\overline{\mathbb{F}}) \to \mathbb{P}^d(W(\overline{\mathbb{F}}))$. In particular, it is infinite.

\item[$\mathrm{(2)}$]
Let $P$ be a fixed height-one prime ideal appearing in $\Char_{R_{W(\overline{\mathbb{F}})}}(M)$. Then one can find an infinite sequence $\{\mathbf{x}_{\widetilde{a}_i}R_{W(\overline{\mathbb{F}})}\}_{i \in \mathbb{N}} \subseteq \mathcal{L}_{W(\overline{\mathbb{F}})}(M)$ 
such that the union:
$$
\bigcup_{i \in \mathbb{N}}\Min_{R_{W(\overline{\mathbb{F}})}}(P+\mathbf{x}_{\widetilde{a}_i}R_{W(\overline{\mathbb{F}})})
$$
is an infinite set.
\end{enumerate}
\end{lemma}

\begin{proof}
$(1)$ By Lemma \ref{lemma7}, we have $\mathbf{x}_{\widetilde{a}} \in \mathcal{L}_{W(\overline{\mathbb{F}})}(M_{\mathrm{null}})$ 
if and only if $\mathbf{x}_{\widetilde{a}} \in \mathcal{L}_{W(\overline{\mathbb{F}})}$ is contained in none of height-two primes $Q'_1,\ldots,Q'_k 
\in \Ass_{R_{W(\overline{\mathbb{F}})}}(M_{\mathrm{null}})$. Let $Q_1,\ldots,Q_{\ell}$ be prime ideals attached to $M$ as in Lemma \ref{lemma6}. 
By the assumption that $\dim R \ge 3$, all these primes together with the set of height-one primes $P_1,\ldots,P_m \in \Supp M$ are strictly contained in the maximal ideal of $R_{W(\overline{\mathbb{F}})}$. So let us find a non-empty Zariski open subset of $\mathbb{P}^d(\overline{\mathbb{F}})$ with the required properties via Lemma \ref{lemma7}. In other words, it suffices to choose $\mathbf{x}_{\widetilde{a}}$ such that 
\begin{equation}
\label{avoidance}
\mathbf{x}_{\widetilde{a}} \notin (\bigcup_{1 \le i \le m} P_i) \cup (\bigcup_{1 \le i \le k} Q'_i) \cup (\bigcup_{1 \le i \le \ell}Q_i).
\end{equation}
Then applying Lemma \ref{Case1Case2} to each prime appearing in $(\ref{avoidance})$, we may find a Zariski open subset $U \subseteq \mathbb{P}^d(\overline{\mathbb{F}})$ such that the condition $(\ref{avoidance})$ holds for $a=(a_0:\cdots:a_d) \in \Sp_{W(\overline{\mathbb{F}})}^{-1}(U)$. Then combining the conclusion of Proposition \ref{prop2}, we complete the proof of (1).

$(2)$ The proof will be completed through inductive steps using $(1)$ as follows. Let $P$ be a fixed height-one prime ideal appearing in $\Char_{R_{W(\overline{\mathbb{F}})}}(M)$. Since $P$ is a non-maximal prime ideal of $R_{W(\overline{\mathbb{F}})}$, we can attach a non-empty Zariski open subset of $\mathbb{P}^d(\overline{\mathbb{F}})$ by Lemma \ref{Case1Case2} and we may find $\mathbf{x}_{\widetilde{a}_0}R_{W(\overline{\mathbb{F}})} \in \mathcal{L}_{W(\overline{\mathbb{F}})}(M)$, which is not contained in $P$. Then this initial choice satisfies our requirement.

Choose $\mathbf{x}_{\widetilde{a}_1}R_{W(\overline{\mathbb{F}})} \in \mathcal{L}_{W(\overline{\mathbb{F}})}(M)$ such that $\mathbf{x}_{\widetilde{a}_1}$ is contained in no prime of the set: 
$$
\Min_{R_{W(\overline{\mathbb{F}})}}(P+\mathbf{x}_{\widetilde{a}_0}R_{W(\overline{\mathbb{F}})})
$$ 
(which is a finite set of primes strictly contained in the maximal ideal of  $R_{W(\overline{\mathbb{F}})}$ due to $\dim R \ge 3$). Next, choose $\mathbf{x}_{\widetilde{a}_2} \in \mathcal{L}_{W(\overline{\mathbb{F}})}(M)$ such that $\mathbf{x}_{\widetilde{a}_2}$ is in no prime ideals contained in the set:
$$
\Min_{R_{W(\overline{\mathbb{F}})}}(P+\mathbf{x}_{\widetilde{a}_0}R_{W(\overline{\mathbb{F}})}) \cup \Min_{R_{W(\overline{\mathbb{F}})}}(P+\mathbf{x}_{\widetilde{a}_1}R_{W(\overline{\mathbb{F}})}). 
$$
By continuing this process, we will eventually obtain a sequence $\mathbf{x}_{\widetilde{a}_0},\mathbf{x}_{\widetilde{a}_1},\mathbf{x}_{\widetilde{a}_2},\ldots$ with the required properties.
\end{proof}

Let $M$ be a finitely generated module over a complete local normal domain $R$ with $W(\mathbb{F})$ its coefficient ring and let $W(\mathbb{F}) \to W(\mathbb{F}')$ be an extension induced by an algebraic extension $\mathbb{F} \to \mathbb{F}'$. We set $M_{W(\mathbb{F}')}:=M \widehat{\otimes}_{W(\mathbb{F})} W(\mathbb{F}')$, which coincides with our previous notation. Let $M$ be a finitely generated torsion $R$-module. Then 
$$
\Char_{R_{W(\mathbb{F}')}}(M_{W(\mathbb{F}')})=\Char_R(M)R_{W(\mathbb{F}')},
$$
which may be verified directly from the definition.

\begin{lemma}
\label{lemma9}
Let $M, N$ be finitely generated torsion $R$-modules. Then we have
$$
\Char_R(M) \subseteq \Char_R(N) \iff \Char_{R_{W(\overline{\mathbb{F}})}}(M_{W(\overline{\mathbb{F}})}) \subseteq \Char_{R_{W(\overline{\mathbb{F}})}}(N_{W(\overline{\mathbb{F}})}).
$$
\end{lemma}

\begin{proof}
The implication $\Rightarrow$ is obvious. So let us prove the other implication. Since $R \to R_{W(\mathbb{F})^{\ur}}$ is ind-\'etale, it suffices to show that
$$
\Char_{R_{W(\mathbb{F})^{\ur}}}(M) \subseteq \Char_{R_{W(\mathbb{F})^{\ur}}}(N) \iff \Char_{R_{W(\overline{\mathbb{F}})}}(M_{W(\overline{\mathbb{F}})}) \subseteq \Char_{R_{W(\overline{\mathbb{F}})}}(N_{W(\overline{\mathbb{F}})}).
$$
Note that $R_{W(\mathbb{F})^{\ur}} \to R_{W(\overline{\mathbb{F}})}$ is faithfully flat, with trivial residue field extension. Then the claim follows from this.
\end{proof}

Now we prove the following theorem.

\begin{theorem}[Control Theorem for Characteristic Ideals]
\label{charid}
With the notation and the hypothesis as in Corollary \ref{cor1}, assume that $M$ and $N$ are finitely generated torsion $R$-modules. Then the following statements are equivalent:

\begin{enumerate}
\item[$\mathrm{(1)}$]
$\Char_R(M) \subseteq \Char_R(N)$.

\item[$\mathrm{(2)}$]
For all but finitely many height-one primes:
$$
\mathbf{x}_{\widetilde{a}}R_{W(\overline{\mathbb{F}})} \in  \mathcal{L}_{W(\overline{\mathbb{F}})}(M_{W(\overline{\mathbb{F}})}) \cap \mathcal{L}_{W(\overline{\mathbb{F}})}(N_{W(\overline{\mathbb{F}})}),
$$
we have
$$
\Char_{R_{W(\mathbb{F}')}/\mathbf{x}_{\widetilde{a}}R_{W(\mathbb{F}')}}(M_{W(\mathbb{F}')}/\mathbf{x}_{\widetilde{a}}M_{W(\mathbb{F}')}) \subseteq \Char_{R_{W(\mathbb{F}')}/\mathbf{x}_{\widetilde{a}}R_{W(\mathbb{F}')}}(N_{W(\mathbb{F}')}/\mathbf{x}_{\widetilde{a}}N_{W(\mathbb{F}')}),
$$
where $\mathbb{F}'$ is any finite field extension of $\mathbb{F}$ depending on ${\widetilde{a}}$ such that $\mathbf{x}_{\widetilde{a}} \in R_{W(\mathbb{F}')}$.

\item[$\mathrm{(3)}$]
For all but finitely many height-one primes:
$$
\mathbf{x}_{\widetilde{a}}R_{W(\overline{\mathbb{F}})} \in  \mathcal{L}_{W(\overline{\mathbb{F}})}(M_{W(\overline{\mathbb{F}})}) \cap \mathcal{L}_{W(\overline{\mathbb{F}})}(N_{W(\overline{\mathbb{F}})}),
$$
we have
$$
\Char_{R_{W(\overline{\mathbb{F}})}/\mathbf{x}_{\widetilde{a}}R_{W(\overline{\mathbb{F}})}}(M_{W(\overline{\mathbb{F}})}/\mathbf{x}_{\widetilde{a}}M_{W(\overline{\mathbb{F}})}) \subseteq \Char_{R_{W(\overline{\mathbb{F}})}/\mathbf{x}_{\widetilde{a}}R_{W(\overline{\mathbb{F}})}}(N_{W(\overline{\mathbb{F}})}/\mathbf{x}_{\widetilde{a}}N_{W(\overline{\mathbb{F}})}).
$$
\end{enumerate} 
\end{theorem}

Note that $\dim R \ge 3$ holds automatically, due to the hypothesis $\depth R \ge 3$.

\begin{proof}
The implications $(1) \Rightarrow (2) \Rightarrow (3)$ are obvious in view of Definition \ref{definition:L}. So it remains to prove $(3) \Rightarrow (1)$. By Lemma \ref{lemma9}, it suffices to show that
$$
\Char_{R_{W(\overline{\mathbb{F}})}}(M_{W(\overline{\mathbb{F}})}) \subseteq \Char_{R_{W(\overline{\mathbb{F}})}}(N_{W(\overline{\mathbb{F}})}).
$$
Take fundamental pseudo-isomorphisms for $M$ and $N$:
$$
M \to \bigoplus_{i} R_{W(\overline{\mathbb{F}})}/P_i^{e_i}~(\mbox{resp}.~N \to \bigoplus_{j} R_{W(\overline{\mathbb{F}})}/Q_j^{f_j})
$$
for a (not necessarily redundant) finite set of height-one primes $\{P_i\}$ (resp. $\{Q_j\}$) of $R_{W(\overline{\mathbb{F}})}$. Put $I_M:=(\prod_i P_i^{e_i})^{\cl}$ and $I_N:=(\prod_j Q_j^{f_j})^{\cl}$ and the condition \textbf{(B)} in Definition \ref{definition:L} allows one to assume that
\begin{equation}
\label{fundamental}
M=\bigoplus_{i} R_{W(\overline{\mathbb{F}})}/P_i^{e_i}~(\mbox{resp}.~N= \bigoplus_{j} R_{W(\overline{\mathbb{F}})}/Q_j^{f_j}).
\end{equation}
To simplify the notation, assume that $\{P_i\}$ (resp. $\{Q_j\}$) is a redundant set of prime ideals and all relevant modules are defined over $R_{W(\overline{\mathbb{F}})}$. We require several steps to complete the proof of the theorem. Let 
$$
\{\mathbf{x}_{\widetilde{a}_i}R_{W(\overline{\mathbb{F}})}\}_{i \in \mathbb{N}} \subseteq \mathcal{L}_{W(\overline{\mathbb{F}})}(M_{W(\overline{\mathbb{F}})}) \cap \mathcal{L}_{W(\overline{\mathbb{F}})}(N_{W(\overline{\mathbb{F}})})
$$ 
be any infinite sequence of distinct primes of $R_{W(\overline{\mathbb{F}})}$ satisfying the condition $(3)$. In particular, we have $\bigcap_{i \in \mathbb{N}} \mathbf{x}_{\widetilde{a}_i}R_{W(\overline{\mathbb{F}})}=0$.

${\bf{Step 1}}$: In this step, we establish $\Supp_{\Ht=1} N \subseteq \Supp_{\Ht=1} M$, where $\Supp_{\Ht=1}(-)$ is the set of height-one primes contained in the support of a module. By assumption, we have
$$
\Big(I_M(R_{W(\overline{\mathbb{F}})}/\mathbf{x}_{\widetilde{a}_i}
R_{W(\overline{\mathbb{F}})})\Big)^{\cl} \subseteq \Big(I_N(R_{W(\overline{\mathbb{F}})}/\mathbf{x}_{\widetilde{a}_i}
R_{W(\overline{\mathbb{F}})})\Big)^{\cl}
$$
for all $i \in \mathbb{N}$. Rewriting this inclusion, we get
$$
\Big((I_M+\mathbf{x}_{\widetilde{a}_i}
R_{W(\overline{\mathbb{F}})})/\mathbf{x}_{\widetilde{a}_i}
R_{W(\overline{\mathbb{F}})} \Big)^{\cl} \subseteq 
\Big((I_N+\mathbf{x}_{\widetilde{a}_i}
R_{W(\overline{\mathbb{F}})})/\mathbf{x}_{\widetilde{a}_i}
R_{W(\overline{\mathbb{F}})} \Big)^{\cl}.
$$
From this description, we deduce the following fact. Fix a height-one prime ideal $Q_k$ from $(\ref{fundamental})$. We may choose the set $\{\mathbf{x}_{\widetilde{a}_i}\}_{i \in \mathbb{N}}$ such that
\begin{equation}
\label{set1}
\bigcup_{i \in \mathbb{N}} \Min_{R_{W(\overline{\mathbb{F}})}}(Q_k+\mathbf{x}_{\widetilde{a}_i}R_{W(\overline{\mathbb{F}})})
\end{equation}
is an infinite set in view of Lemma \ref{lemma8}. For every fixed $i \in \mathbb{N}$, we have
$$
I_M \subseteq \fp_i,
$$ 
where $\fp_i \in \Min_{R_{W(\overline{\mathbb{F}})}}(Q_k+\mathbf{x}_{\widetilde{a}_i}R_{W(\overline{\mathbb{F}})})$ is chosen to be an arbitrary fixed element. 

On the other hand, since $(\ref{set1})$ is infinite, we may find an infinite subset $\{\fp_i\}_{i \in \mathbb{N}}$ of $(\ref{set1})$ and we fix it once and for all. Since $R_{W(\overline{\mathbb{F}})}/Q_k$ is an integral domain and $\{\fp_i\}_{i \in \mathbb{N}}$ is an infinite set of height-two primes containing $Q_k$, it follows from Lemma \ref{Noetherian} that
$$
Q_k=\bigcap_{i \in \mathbb{N}} \fp_i.
$$
Since $I_M \subseteq \fp_i$ for all $i \in \mathbb{N}$, we have $I_M \subseteq Q_k$. Since $Q_k$ is arbitrary, $I_M \subseteq (\prod_j Q_j)^{\cl}$, or equivalently, $\Supp_{\Ht=1} N \subseteq \Supp_{\Ht=1} M$.

${\bf{Step 2}}$: In this step, we deal with multiplicities of divisors in the characteristic ideal and we complete this step by induction on the number of divisors appearing in $I_M$.

First, assume $I_M=(P^e)^{\cl}$ for $e \ge 1$. Then we have $\Supp_{\Ht=1} N=\varnothing$ or $\{P\}$ because of ${\bf{Step 1}}$. If $\Supp_{\Ht=1} N=\varnothing$, there is nothing to prove. So assume $\Supp_{\Ht=1} N=\{P\}$. Both $M$ and $N$ are assumed to be fundamental torsion $R_{W(\overline{\mathbb{F}})}$-modules, thus $M[\mathbf{x}_{\widetilde{a}_i}]$ and $N[\mathbf{x}_{\widetilde{a}_i}]$ are trivial modules and Proposition \ref{prop3} yields that $\ell_{(R_{W(\overline{\mathbb{F}})})_P}(M_P) \ge \ell_{(R_{W(\overline{\mathbb{F}})})_P}(N_P)$.

In the general case, we prove by contradiction and thus, assume that $I_M \nsubseteq I_N$. Then this implies that we have $e_k<f_k$ for some $k$, where $e_k,f_k$ are coming from $(\ref{fundamental})$. Put
\begin{equation}
\label{reflexiveideal}
\widetilde{I}_M:=(P_k^{-e_k} \cdot I_M)^{\cl}~(\mbox{resp}.~\widetilde{I}_N:=(P_k^{-e_k} \cdot I_N)^{\cl}),
\end{equation}
which are both integral reflexive ideals. There are short exact sequences:
$$
0 \to \widetilde{I}_M/I_M \to R_{W(\overline{\mathbb{F}})}/I_M \to R_{W(\overline{\mathbb{F}})}/\widetilde{I}_M \to 0
$$
and
$$
0 \to \widetilde{I}_N/I_N \to R_{W(\overline{\mathbb{F}})}/I_N \to R_{W(\overline{\mathbb{F}})}/\widetilde{I}_N \to 0
$$
and it is clear that
$$
\Char_{R_{W(\overline{\mathbb{F}})}}(\widetilde{I}_M/I_M)=\Char_{R_{W(\overline{\mathbb{F}})}}(\widetilde{I}_N/I_N)=(P_k^{e_k})^{\cl},
$$
which induces the following short exact sequences by the snake lemma:
$$
0 \to \widetilde{I}_M/(I_M,\mathbf{x}_{\widetilde{a}_i}\widetilde{I}_M) \to R_{W(\overline{\mathbb{F}})}/(I_M,\mathbf{x}_{\widetilde{a}_i}) \to R_{W(\overline{\mathbb{F}})}/(\widetilde{I}_M,\mathbf{x}_{\widetilde{a}_i}) \to 0
$$
and
$$
0 \to \widetilde{I}_N/(I_N,\mathbf{x}_{\widetilde{a}_i}\widetilde{I}_N) \to R_{W(\overline{\mathbb{F}})}/(I_N,\mathbf{x}_{\widetilde{a}_i}) \to R_{W(\overline{\mathbb{F}})}/(\widetilde{I}_N,\mathbf{x}_{\widetilde{a}_i}) \to 0.
$$
Taking characteristic ideals, we get from the condition $(3)$ that
$$
\Char_{R_{W(\overline{\mathbb{F}})}/\mathbf{x}_{\widetilde{a}_i}R_{W(\overline{\mathbb{F}})}}(R_{W(\overline{\mathbb{F}})}/(\widetilde{I}_M,\mathbf{x}_{\widetilde{a}_i})) \subseteq \Char_{R_{W(\overline{\mathbb{F}})}/\mathbf{x}_{\widetilde{a}_i}R_{W(\overline{\mathbb{F}})}}(R_{W(\overline{\mathbb{F}})}/(\widetilde{I}_N,\mathbf{x}_{\widetilde{a}_i})).
$$
Since the number of primes ideals in $\Ass_{R_{W(\overline{\mathbb{F}})}}(R_{W(\overline{\mathbb{F}})}/\widetilde{I}_M)$ is 
exactly one less than that of components of prime ideals in $I_M$, the induction hypothesis on $\widetilde{I}_M$ yields that
$$
\Char_{R_{W(\overline{\mathbb{F}})}}(R_{W(\overline{\mathbb{F}})}/\widetilde{I}_M) \subseteq \Char_{R_{W(\overline{\mathbb{F}})}}(R_{W(\overline{\mathbb{F}})}/\widetilde{I}_N).
$$
However, we deduce from these observations and $(\ref{reflexiveideal})$ that $I_M \subseteq I_N$, which is a contradiction to our assumption $I_M \nsubseteq I_N$. Hence, we obtain $I_M \subseteq I_N$, as desired.
\end{proof}

\begin{remark}
It is worth pointing out that Theorem \ref{charid} holds for complete local normal rings of mixed characteristic with arbitrary infinite perfect residue field as well. More precisely, it can be proven that $\Char_R(M) \subseteq \Char_R(N) \iff \Char_{R/xR}(M/xM) \subseteq \Char_{R/xR}(N/xN)$ for sufficiently many $x \in R$.
\end{remark}

In this article, we presented an application of the local Bertini theorem to characteristic ideals. However, we believe that the main theorem has more interesting applications such as the study of the restriction map on divisor class (Chow) groups.

\section{Appendix}\label{appendix}

In this appendix, we study the relationship between Fitting ideals and characteristic ideals. For Fitting ideals, we refer the reader to Northcott's book \cite{Nor}, but we review the basic part of the theory. For reflexive sheaves on normal schemes, we refer te reader to \cite{Har}. Throughout, we assume that $R$ is a Noetherian ring and $M$ is a finitely generated $R$-module.

\begin{definition}[Fitting ideal]
Let the notation be as above and assume that
$$
F_1 \to F_0 \to M \to 0
$$
is a finite free resolution of the $R$-module $M$, where the mapping $F_1 \to F_0$ is defined via a $m \times n$-matrix $X$ with $\rank(F_1)=n$ and $\rank(F_0)=m$. Then $\Fitt_R(M)$ is defined as an ideal of $R$ generated by all $m$-minors of $X$.
\end{definition}

The Fitting ideal does not depend on the choice of a free resolution and it enjoys the following properties.

\begin{proposition}
\label{prop4}
Let $M$ be a finitely generated module over a Noetherian ring $R$. Then we have the following properties.

\begin{enumerate}
\item[$\mathrm{(1)}$]
Let $I \subseteq R$ be an ideal. Then $\Fitt_R(R/I)=I$.

\item[$\mathrm{(2)}$]
Let $S$ be any Noetherian $R$-algebra. Then $\Fitt_S(M \otimes _R S)=\Fitt_R(M)S$.

\item[$\mathrm{(3)}$]
Let $\Ann_R(M)$ be the annihilator of the $R$-module $M$. Then
$$
\Fitt_R(M) \subseteq \Ann_R(M).
$$

\item[$\mathrm{(4)}$]
If $0 \to L \to M \to N \to 0$ is a short exact sequence of $R$-modules, then
$$
\Fitt_R(L) \cdot \Fitt_R(N) \subseteq \Fitt_R(M).
$$

\item[$\mathrm{(5)}$]
Assume that $R$ is a discrete valuation ring with its uniformizing parameter $b$ and $M$ is a torsion $R$-module. Then $\Fitt_R(M)=(b)^{\ell_R(M)}$.
\end{enumerate}
\end{proposition}

\begin{proof}
These facts are all well known. For $(5)$, it simply follows from the elementary divisors of modules over a principal ideal domain.
\end{proof}

For a Noetherian domain $R$ and an $R$-module $M$, let $M^*:=\Hom_R(M,R)$, the dual of $M$. We say that $M^{\cl}:=(M^*)^*$ is the \textit{reflexive closure} of $M$. Then we have the following lemma.

\begin{lemma}
\label{10}
Let $R$ be a Noetherian domain and let $I$ be a fractional ideal of $R$. Then the reflexive closure $I^{\cl}$ is naturally regarded as a fractional ideal of $R$.\end{lemma}

\begin{proof}
By assumption, there exists $\alpha \in R$ such that $I \simeq \alpha \cdot I \subseteq R$. Let $J:=\alpha \cdot I$, an ideal of $R$. The short exact sequence $0 \to J \to R \to R/J \to 0$ induces a short exact sequence
$$
0=\Hom_R(R/J,R) \to R \to \Hom_R(J,R) \to N \to 0,
$$
with $N \subseteq \Ext^1_R(R/J,R)$ cokernel of $R \to \Hom_R(J,R)$. Then applying $\Hom_R(-,R)$ twice, we get an exact sequence
$$
0=\Hom_R(N,R) \to J^{\cl} \to R,
$$
because $J \cdot N=0$. This implies that $J^{\cl}$ is an ideal of $R$. Then $I^{\cl}=\alpha^{-1} \cdot J^{\cl}$ is a fractional ideal.
\end{proof}

Note that even when $I$ and $J$ are reflexive, the product $I \cdot J$ need not be reflexive. Any principal ideal is reflexive. Let $I$ be an ideal of a normal domain $R$. Then we recall the following fact:
$$
I^{\cl}=\bigcap_{P}I_P,
$$
where $P$ ranges over all height-one primes of $R$. The natural inclusion $I \to I^{\cl}$ is a pseudo-isomorphism, since $I_P \to (I^{\cl})_P=(I_P)^{\cl}$ for every height-one prime $P \subseteq R$ and every ideal in a discrete valuation ring is principal. The following lemma explains the naturality of reflexive ideals and gives a way to investigate the inclusion relation between characteristic ideals.

\begin{lemma}
\label{lemma11}
Let $R$ be a Noetherian normal domain and let $I$ and $J$ be reflexive ideals. Then $I \subseteq J$ if and only if $v_P(I) \ge v_P(J)$ for a  valuation $v$ attached to every height-one prime $P$ of $R$. In particular, the only reflexive integral ideal containing a prime ideal of $R$ properly is $R$ itself.
\end{lemma}

If $R$ is only assumed to be Cohen-Macaulay, a similar result holds for invertible modules (\cite{Flach}, Lemma 5.3). We defined characteristic ideals as reflexive ideals and this is natural from the viewpoint of Iwasawa's main conjecture, because the most interesting arithmetic information may be captured at height-one primes. For finitely generated torsion $R$-modules $M, N$, it follows from the above lemma that $\Char_R(M) \subseteq \Char_R(N)$ if and only if $\Char_R(M)_P \subseteq \Char_R(N)_P$ for every height-one prime $P \in \Supp(M) \cup \Supp(N)$.

\begin{example}
Suppose $I$ is reflexive and let $\fa \subseteq R$ be such that $R/\fa$ is a normal domain. Then $I(R/\fa)$ need not be reflexive. For a general ideal $I \subseteq R$, it can happen that
$$
(I(R/\fa))^{\cl} \ne (I^{\cl}(R/\fa))^{\cl}.
$$
Here, $I^{\cl}$ is the reflexive closure with respect to $R$ and $(I(R/\fa))^{\cl}$ is the reflexive closure with respect to $R/\fa$. Let us take a look at the following simple example. Take $R=\mathbb{Z}_p[[x,y]]$, $I=(x,y)$, and $\fa=(x)$. Then $(I^{\cl}(R/\fa))^{\cl}=R/xR$, since there is no height-one prime of $R$ containing $I$. But $(I(R/\fa))^{\cl}=y(R/xR)$. 
\end{example}

\begin{proposition}
\label{prop5}
Let $R$ be a Noetherian normal domain. Then the following hold.

\begin{enumerate}
\item[$\mathrm{(1)}$] 
Let $M$ be a finitely generated torsion $R$-module. Then we have
$$
\Char_R(M)=\big(\prod_{\Ht \fp=1} \fp^{\ell_{R_{\fp}}(M_{\fp})}\big)^{\cl}=\Fitt_R(M)^{\cl}.
$$
In particular, $\Fitt_R(M) \subseteq \Char_R(M)$ and if $R$ is a UFD, then
$$
\Fitt_R(M) \subseteq \prod_{\Ht \fp=1} \fp^{\ell_{R_{\fp}}(M_{\fp})}=\Char_R(M).
$$

\item[$\mathrm{(2)}$]
Let $0 \to L \to M \to N \to 0$ be a short exact sequence of finitely generated $R$-modules. Then
$$
\Char_R(M)=(\Char_R(L) \cdot \Char_R(N))^{\cl}.
$$
\end{enumerate}
\end{proposition}

\begin{proof}
$(1)$ 
Since the characteristic ideal is reflexive, the first equality follows by taking localization at all height-one primes of $R$. The second equality follows from the fact
$$
\Fitt_{R_P}(M_P)=(PR_P)^{\ell_{R_P}(M_P)}
$$
for any height-one prime ideal $P \subseteq R$. The second assertion is due to the fact that a height-one prime in a UFD is principal.

$(2)$ This is clear, since the length is additive with respect to short exact sequences. The equality continues to hold true without torsion property of modules.
\end{proof}

\begin{example}
The ordinary power of a height-one prime in a normal domain is not necessarily reflexive. Here is an example. Let $R=\mathbb{Z}_p[[x^2, xy, y^2]]$ and let $\fp= (x^2, xy)$. Then $R$ is a normal domain and $\Ht(\fp)=1$. Then $\fp^2=(x^4,x^3y,x^2y^2)$ and $\fp^2 \ne \fp^{(2)}$. In fact, $x^2 \notin \fp^2$, but $\fp^{(2)}=(x^2)$. 

Now let $\fq=(xy,y^2)$ and $M=R/(\fp \cap \fq)$, which is torsion over $R$. Then one verifies that 
$$
\Fitt_R(M)=\fp\cap \fq \nsubseteq \prod_{\Ht \fp=1} \fp^{\ell_{R_{\fp}}(M_{\fp})}=\fp\fq,
$$
which tells us that Proposition \ref{prop5} (1) is the most optimal.
\end{example}

\begin{acknowledgement}
The authors thank the anonymous referees for reading the manuscript carefully. The second-named author thanks Prof. Trivedi for explaining the proof of a lemma in \cite{Fle}.
\end{acknowledgement}

\end{document}